\documentclass[a4paper,11pt]{article}

\usepackage{graphicx}
\usepackage{tikz} 
\usepackage{tkz-euclide}
\usepackage{svg}
\usepackage{amsmath}
\usepackage{amsthm}
\usepackage{amssymb}
\usepackage{fancyhdr}
\usepackage{enumerate}
\usepackage[breaklinks]{hyperref}
\usepackage{lscape}
\usepackage{multirow}
\usepackage{relsize}
\usepackage{setspace}
\usepackage[comma]{natbib}
\usepackage{fullpage}
\usepackage{caption}
\usepackage{authblk}
\usepackage{totcount}
\usepackage{array}
\usepackage{longtable}
\usepackage{xcolor}
\usepackage[titletoc]{appendix}
\usepackage[top    = 0.75in,
bottom = 1.22in,
left   = 1in,
right  = 1in]{geometry}
\usepackage{booktabs}

\renewcommand{\vec}[1]{\mathbf{#1}}
\newcommand{\black}[1]{{\color{black}#1}}

\newcolumntype{L}[1]{>{\raggedright\arraybackslash}p{#1}}

\newcolumntype{C}[1]{>{\centering\arraybackslash}m{#1}}

\newtotcounter{citnum} 
\def\oldbibitem{} \let\oldbibitem=\bibitem
\def\bibitem{\stepcounter{citnum}\oldbibitem}

\newtheorem{definition}{Definition}[section]
\newtheorem{theorem}[definition]{Theorem}
\newtheorem{lemma}[definition]{Lemma}
\newtheorem{corollary}[definition]{Corollary}

\newtheorem{remark}[definition]{Remark}

\begin{document}
\onehalfspacing

\title{Stochastic dynamic job scheduling with interruptible setup and processing times\black{: An approach based on queueing control}}


\author{\large Dongnuan Tian\\ \vspace{-4mm} \footnotesize{Department of Management Science, Lancaster University Management School, Lancaster University, Lancaster LA1 4YW, United Kingdom.  Email: d.tian2@lancaster.ac.uk.}\\\text{ }\\ \large Rob Shone\footnote{Corresponding author}\\ \footnotesize{Department of Management Science, Lancaster University Management School, Lancaster University, Lancaster LA1 4YW, United Kingdom.  Email: r.shone@lancaster.ac.uk.}\\}
\date{ }

\pagestyle{fancy}
\fancyhf{}
\rhead{\thepage}
\lhead{\footnotesize \textbf{Tian and Shone: }\textit{Stochastic dynamic job scheduling with interruptible setup and processing times}}



\normalsize
\maketitle

\vspace{-6mm}
\hrule
\text{ }\\

\noindent \Large \textbf{Abstract}\\
\normalsize

\noindent We consider a stochastic, dynamic job scheduling problem, \black{formulated as a queueing control problem,} in which a single server processes jobs of different types that arrive according to independent Poisson processes. The problem is defined on a network, with jobs arriving at designated demand points and waiting in queues to be processed by the server, which travels around the network dynamically and is able to change its course at any time. In the context of machine scheduling, this enables us to consider sequence-dependent, interruptible setup and processing times, with the network structure encoding the amounts of effort needed to switch between different tasks. We formulate the problem as a Markov decision process in which the objective is to minimize long-run average holding costs and prove the existence of a stationary policy under which the system is stable, subject to a condition on the workload of the system. We then propose a class of index-based heuristic policies, show that these possess intuitively appealing structural properties and suggest how to modify these heuristics to ensure scalability to larger problem sizes. Results from extensive numerical experiments are presented in order to show that our heuristic policies perform well against suitable benchmarks. \\

\noindent \textit{Keywords: }Job scheduling; dynamic programming; index heuristics\\

\hrule


\section{Introduction}
\label{sec:Introduction}

\text{ }\indent Job scheduling and server scheduling problems have been very widely studied in operations research. 
Stochastic, dynamic versions of such problems typically draw upon both scheduling and queueing theory (\cite{leung2004handbook, pinedo2016scheduling, blazewicz2019scheduling}). 
A classical problem formulation involves a system of $N$ parallel queues, each with its own arrival process, and a single server with the ability to switch dynamically between queues. Jobs at queue $i$ \black{arrive at a rate of $\lambda_i$ and can be processed at a rate of $\mu_i$, $i=1,...,N$. A holding cost $c_i$ is incurred for each unit of time that a type $i$ job spends in the system.} The server is able to observe queue lengths continuously and can provide service to one queue at a time. This type of formulation induces a well-known $c\mu$-rule, whereby the queues are ranked in descending order of the product $c_i\mu_i$, and the server always selects a job from the queue with the largest $c_i\mu_i$ value among the non-empty queues. 
The optimality of the $c\mu$-rule in various queueing settings has been well-documented (\cite{smith1956various,baras1985two,buyukkoc1985cmu,vanmieghem1995}).

In recent decades, researchers have sought to extend the applicability of the $c\mu$-rule to address more complex scheduling problems. \cite{mandelbaum2004scheduling} introduced a generalized version of the rule for systems with convex delay cost structures, reflecting more realistic scenarios where penalties for delays increase at an accelerating rate. \cite{atar2010cmu} modified the classical $c\mu$-rule to incorporate abandonment effects, resulting in the $c\mu/\theta$-rule. This rule considers both service rates and abandonment rates, prioritizing queues based on the ratio of the cost-weighted service rate to the abandonment rate. 
\cite{veatch2015cmu} proposed a $c\mu$-rule for a parallel flexible server system with a two-tier structure. In this model, the first tier consists of job classes that are assigned to specific servers, while the second tier includes a class of jobs that can be served by any available server. 
\black{In more recent years, applications of the $c\mu$-rule have continued to attract a lot of attention.} \cite{lee2021scheduling} established a learning-based variant of the rule, where the algorithm learns the job parameters over time and adapts its scheduling decisions to minimize cumulative holding costs, even when the statistical parameters of the jobs are initially unknown. \black{\cite{cohen2022asymptotic} also considered uncertain model parameters, and showed that the $c\mu$-rule attains a type of asymptotic optimality. \cite{ozkan2022cmu} investigated the performance of the $c\mu$-rule in systems with non-parallel configurations, such as tandem queues.} 

Job scheduling and server scheduling problems have broad applications in several real-world domains including cloud computing, communications networks and manufacturing (\cite{shaw2014survey,zhang2019resource,xu2021scheduling}). In this paper we consider a stochastic, dynamic job scheduling problem formulated on a network of nodes and edges, in which certain nodes are designated as `demand points' and act as entry points for jobs of specific types. There is a single server, which may represent a machine (for example) that needs to process jobs one at a time. In this context, the time needed for the server to move from one demand point to another represents the setup time required for the machine to switch from processing one type of job to another. The design of the network may be seen as a way to encode the amounts of effort needed for the server (or machine) to switch between jobs of different types.

\begin{figure}[hbtp]
    \begin{center}
        \includeinkscape[scale = 0.8]{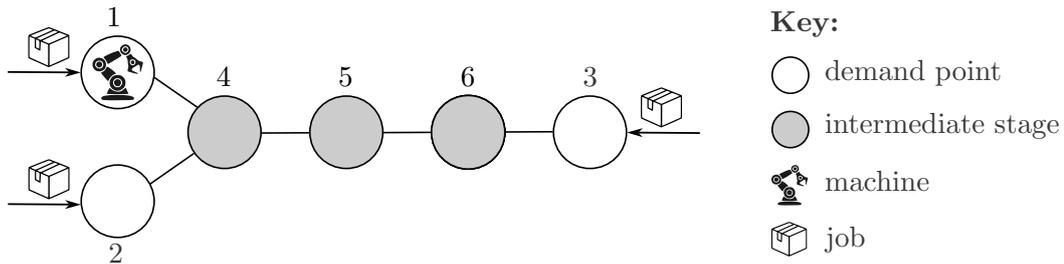}
    \end{center}
    \caption{A network with 3 demand points, 3 intermediate stages and a single server. White-colored nodes represent demand points at which new jobs arrive, and gray-colored nodes represent intermediate stages.}
    \label{fig1}
\end{figure}

As an example of our network-based approach, consider the situation shown in Figure \ref{fig1}. The white-colored nodes are demand points, at which new jobs arrive according to independent Poisson processes, and the gray-colored nodes are `intermediate stages', which must be traversed in order to move from either of the left-hand demand points to the right-hand demand point or vice versa. A single server occupies one node in the network at any given time, and can either remain at its current node or attempt to move to an adjacent node. In order to process jobs of type $i$ (for $i\in\{1,2,3\}$), it must be located at node $i$. The times required to process jobs and switch between adjacent nodes are random, and jobs waiting to be processed incur linear holding costs. We provide a detailed problem formulation in Section \ref{sec:Problem Formulation}.

An important feature of our network formulation is the fact that the server is always free to choose a new action (i.e. a new direction to move), regardless of actions chosen at previous time points. This implies that it can interrupt the processing of a particular job, or the transition between different job types, in order to follow a different course. For example, suppose the server is located at node 1 in Figure \ref{fig1}, but wishes to move to node 3 so that it can process jobs of type 3. To accomplish this, it must pass through the intermediate nodes 4, 5 and 6. We assume that the amount of time needed to switch between any two adjacent nodes is randomly distributed (further details are provided in Section \ref{sec:Problem Formulation}). However, suppose that after arriving at node 5, it decides to change course and move towards node 2 instead. (This may happen, for example, if new jobs have recently arrived at node 2.) In the job scheduling context, we would say that the setup time needed to prepare to process jobs of type 3 is interrupted so that the server can prepare to process jobs of type 2 instead. Moreover, our network formulation not only allows setup times to be interrupted, but also allows them to depend on any partial progress made on interrupted setup attempts. For example, in Figure \ref{fig1}, job types 1 and 2 may be regarded as having similar setup requirements, as they are located close to each other in the network. If the server begins at node 3 and initially moves towards node 1, but then decides to move to node 2 instead, then any progress made in setting up jobs of type 1 (i.e. moving towards node 1) is helpful in reducing the amount of time needed for setting up jobs of type 2. On the other hand, if the server begins at node 1 and initially moves towards node 3, but then decides to move to node 2 instead, then any progress made in setting up jobs of type 3 actually delays the setting up of type 2 jobs. \black{Accurate modeling of complex setup time requirements (including the effects of interrupted setups) is important in settings where new jobs arrive frequently and randomly, requiring the server to allocate its time efficiently.} 

While our network formulation is intuitively quite simple, it allows consideration of features that are not particularly common in job scheduling problem formulations, such as sequence-dependent, interruptible setup times that are affected by partial setup progress. These kinds of generalizations can be particularly relevant in production and manufacturing applications, where the reallocation of a server from one task to another can involve complicated reconfigurations or the assembly of specialized tools and equipment (\cite{allahverdi1999setup,gholami2009setup}). We also note that although job scheduling applications are a primary motivation for our study, our formulation is sufficiently general to allow potential applications to other common operations research problems. For example, one could consider a dynamic vehicle routing problem in which a vehicle provides service to customers in geographically distinct locations \black{(\cite{ulmer2020vrp})}, or a security problem in which a defensive agent responds to threats that appear in a computer network \black{(\cite{hunt2024cyber})}. 

The stochastic, dynamic nature of our problem implies that a Markov decision process (MDP) formulation is appropriate, although (as in many other problems formulated as MDPs) it is not possible to compute optimal solutions unless the scale of the problem is very small; thus, we must consider heuristic approaches. Several previous studies have also used MDP formulations in job scheduling contexts, although their assumptions are usually quite different from ours. \cite{8248100} considered a problem in which jobs with different processing requirements arrive randomly and wait in queues to be processed by machines, and at each decision epoch a job must be chosen to be processed next. The machines are relatively `static' in this setting, rather than being assigned dynamically to different tasks. \cite{luo2020dynamic} \black{(see also \cite{lei2023jobshop,Zhao2023jobshop})} formulated a dynamic flexible job shop scheduling problem with new job insertions and used a Q-learning approach, but their objective is quite different from ours as it is based on minimizing the total tardiness of a given set of jobs with fixed processing times. \cite{elsayed2022deep} also considered the processing of a given set of jobs on multiple machines and used a novel graph-based MDP formulation, but their study does not consider the stochastic factors present in our model. \cite{fan2012supply} considered a single machine scheduling problem (SMSP) with transportation costs and used an MDP algorithm to minimize the overall cost of processing jobs in a sequence and transporting the finished products to a single customer; however, their problem (unlike ours) is based on a finite set of jobs. \cite{yang2022dynamic} formulated an SMSP in which the machine's state is subject to uncertainty, resulting in job time parameters being expressed in probabilistic terms. Their problem is also of a finite-time nature, with an objective based on minimizing the makespan. \black{In summary, our paper addresses an important research gap by making a connection between job scheduling problems with random arrivals and infinite-horizon queueing control problems.}.

A particularly relevant paper to ours is \cite{duenyas1996heuristic}, which considers the problem of dynamically allocating a single server to process jobs of different types that arrive according to independent Poisson processes and wait in parallel queues. As in our problem, the authors consider random processing times and switching times (between different queues), and the time horizon of interest is infinite. However, they do not allow processing times or switching times to be interrupted. Furthermore, they do not use a network formulation, and the time required to switch from one queue to another depends only on the destination queue; thus, there are no sequence-dependent setup times. Despite these differences, their approach of developing index heuristics is similar to the approach that we pursue in this study, and their algorithms use certain steps and conditions that we aim to adapt and generalize for use in our problem (full details can be found in Section \ref{sec:Index heuristic}). We also note that the study of \cite{duenyas1996heuristic} makes use of certain results from the classical literature on `polling systems', in which a single server visits a set of queues according to a predefined sequence (\cite{boxma1987pseudo,browne1989dynamic,altman1992stability,yechiali1993analysis}). In Section \ref{sec:Problem Formulation}, we also draw upon this body of research to show that our system possesses an important stability property.

As an alternative to index heuristics, one can also consider reinforcement learning (RL) approaches, and several recent studies have applied RL to machine scheduling problems (\cite{li2020reinforcement,wang2021review,kayhan2023reinforcement,li2024reinforcement}). Although RL methods are undoubtedly powerful, they also have some drawbacks in comparison to alternative approaches (\cite{dulac2019challenges}). Solutions given by RL algorithms tend to be less interpretable than those given by simpler heuristics, and therefore less appealing to system operators. Additionally, there is the issue of online replanning, as discussed in \cite{bertsekas2019rl}. If the parameters of the problem (e.g. job arrival rates or setup time distributions) change suddenly and unexpectedly, RL methods may struggle to adapt `on-the-fly' as they require expensive offline training in order to be able to discover strong-performing policies, whereas index heuristics can quickly adapt to the new parameters of the problem. Thus, whilst we acknowledge the potential of RL methods, these are not within the scope of our current study.

The main contributions of this paper are as follows:\\

\vspace{-5mm}
\begin{itemize}
\item \black{We provide a novel formulation of a stochastic, dynamic job scheduling problem in which the setup time distributions for different types of jobs are encoded by a network structure, allowing them to be both sequence-dependent and interruptible.}
\item We prove that there always exists a decision-making policy under which the system is stable, subject to a condition on the total workload of the system.
\item We propose a class of index-heuristics, referred to as $K$-stop heuristics, which are suitable for our network-based problem as they make decisions by taking the topological structure of the network into account.
\item We prove that our heuristics possess a pathwise consistency property which ensures that the server always proceeds to a demand point in finite time. \black{Moreover, we show that these heuristics attain both system stability and optimality under certain conditions.}
\item We propose a further class of heuristics, known as ($K$ from $L$)-stop heuristics, that scale much more readily to large network sizes than the $K$-stop heuristics.
\item We present extensive results from numerical experiments in order to compare the relative performances of our heuristics, comparing them (where possible) to optimal values given by dynamic programming and also to the performance of the unmodified heuristic policy in \cite{duenyas1996heuristic}.\\
\end{itemize}
\vspace{-5mm}

The rest of the paper is organized as follows. In Section \ref{sec:Problem Formulation} we formulate our stochastic, dynamic job scheduling problem as an MDP and prove an important stability property. In Section \ref{sec:Index heuristic} we derive index heuristics and prove some useful properties of these heuristics, such as pathwise consistency. In Section \ref{sec:numerical} we present the results of our numerical study. Finally, our concluding remarks can be found in Section \ref{sec:conclusions}.\\

\section{Problem Formulation}
\label{sec:Problem Formulation}
\text{ }\indent The problem is defined on a connected graph, referred to as a \emph{network}. Let $V$ and $E$ denote the sets of nodes and edges respectively, and let $D\subseteq V$ be a subset of nodes referred to as \emph{demand points} (the white nodes in Figure \ref{fig1}). Let $N:=V\setminus D$ denote the other nodes (the gray nodes in Figure \ref{fig1}), referred to as \emph{intermediate stages}. We use $d=|D|$ and $n=|N|$ to denote the numbers of demand points and intermediate stages, respectively. We will also assume that the demand points are numbered $1,2,...,d$ and the intermediate stages are numbered $d+1,d+2,...,d+n$. Jobs arrive at demand point $i\in D$ according to a Poisson process with intensity rate $\lambda_i>0$, referred to as an \emph{arrival rate}, and wait in a \black{first-come-first-served} queue until they are processed. Arrivals at different demand points are assumed to occur independently of each other. Additionally, a linear holding cost $c_i>0$ per unit time is incurred \black{for each job waiting to be processed} at node $i \in D$. 

Jobs are processed by a single \emph{server} (or \emph{machine}, \emph{resource} etc.) which can move around the network and, at any given time, is located at a single node in $V$. At any point in time, the server can either remain at its current node or make an attempt to move to an adjacent node. In the latter case, the server must also decide which node to move to. Thus, if the server's current node is adjacent to $k$ other nodes then there are $k+1$ possible decisions for the server. If the server decides to remain at a node $i\in D$ (i.e. at a demand point) and this node has a number of jobs $x_i>0$ waiting to be processed, 
then the jobs are processed at an instantaneous rate $\mu_i>0$, where $\mu_i$ is the \emph{processing rate} for demand point $i$. If the server remains at a demand point $i$ with no jobs present (i.e. $x_i=0$), then the server is said to be \emph{idle}. Likewise, idleness can also occur when the server chooses to remain at some intermediate stage $j\in N$. On the other hand, if the server chooses to move (or `switch') to an adjacent node, then the switch occurs at an instantaneous rate $\tau>0$, referred to as a \emph{switching rate}. Switching and processing times are assumed to be independent of the arrival processes for demand points in $D$.

The assumption of instantaneous processing and switching rates implies that processing and switching times are exponentially distributed, but they are also \emph{interruptible} because the server can change its decision at any point in time. For example, it can choose to remain at a non-empty demand point but then choose to switch before any further jobs have been processed. \black{Alternatively, it could choose to switch from one node $i$ to another node $j$, but then change direction and attempt to switch to a different node $k$ before the switch to $j$ is complete. We note that, 
although the switching time between two adjacent nodes has the memoryless property in our model, the time to switch from one demand point to another (which, in general, requires passing through a sequence of intermediate stages) is instead distributed as a sum of i.i.d. exponential switching times, implying that it has an Erlang distribution. Throughout this paper we assume exponential processing times and switching times (between adjacent nodes), with the latter assumption implying Erlang-distributed setup times (between demand points). However, later in this section we also comment on how the formulation could be extended to the more general case of phase-type processing and switching times.}

Under the above assumptions, the system can be formulated as a continuous-time Markov decision process (MDP). The state space can be written as
$$S:=\left\{(v,(x_1,...,x_d))\;|\;v\in V,\;x_i\geq 0 \text{ for }i\in D\right\},$$
where $v$ is the node currently occupied by the server \black{and $x_i$ is the number of jobs waiting to be processed (including any job currently being processed) at demand point $i\in D$}. We will use vectors such as  $\vec{x}$ and $\vec{y}$ to represent generic states in $S$ and use $v(\vec{x})$ to denote the server's location under state $\vec{x}\in S$. In order to simplify notation we will sometimes write $v$ instead of $v(\vec{x})$ if the state associated with $v$ is clear. Considering that the server cannot be processing jobs whilst also switching and also cannot serve more than one demand point at a time, the total of the transition rates under any state is at most $\sum_{i\in D}\lambda_i+\max\{\mu_1,...,\mu_d,\tau\}$ and we can therefore use the technique of uniformization (\cite{serfozo1979equivalence}) to transform the system into a discrete-time counterpart that evolves in time steps of size
\begin{equation}\Delta:=\left(\sum_{i\in D}\lambda_i+\max\{\mu_1,...,\mu_d,\tau\}\right)^{-1}.\label{Delta_eqn}\end{equation}

Let $R(\vec{x})$ be the set of nodes adjacent to node $v(\vec{x})$ \black{(not including $v(\vec{x})$ itself)} and $A_\vec{x}=\{v(\vec{x})\}\cup R(\vec{x})$ be the action set available under state $\vec{x}\in S$ at a particular time step. We can interpret action $a\in A_\vec{x}$ as the node that the server tries to move to next, with $a=v(\vec{x})$ indicating that the server remains where it is. Then, for any pair of states $\vec{x},\vec{y}\in S$ with $\vec{x}\neq \vec{y}$, the transition probability of moving from $\vec{x}:=(v(\vec{x}),(x_1,...,x_d))$ to $\vec{y}:=(v(\vec{y}),(y_1,...,y_d))$ following the choice of action $a \in A_{\vec{x}}$ can be expressed as
\begin{equation}p_{\vec{x},\vec{y}}(a):=\begin{cases}\lambda_i\Delta,&\text{ if }y_i=x_i+1\text{ and }y_j=x_j\text{ for }j\neq i,\\[8pt]
	\mu_i\Delta,&\text{ if }y_i=x_i-1,\;y_j=x_j\text{ for }j\neq i\text{ and }a=v(\vec{x})=i,\\[8pt]
	\tau\Delta,&\text{ if }y_j=x_j\text{ for all }j\in D,\;v(\vec{y})\in R(\vec{x})\text{ and }a=v(\vec{y}),\\[8pt]
	0,&\text{ otherwise,}\end{cases}\label{trans_prob_eqn}\end{equation}
with $p_{\vec{x},\vec{x}}(a)=1-\sum_{\vec{y} \neq \vec{x}}p_{\vec{x},\vec{y}}(a)$ being the probability of remaining at the same state. \black{Since the units of time are arbitrary, we will assume $\Delta=1$ without loss of generality. This implies that parameters such as $\lambda_i$, $\mu_i$ and $\tau$ can be interpreted as probabilities in our uniformized MDP, rather than transition rates.}

At each time step, the single-step cost $f(\vec{x})$ is calculated by summing the holding costs of jobs that are still waiting to be processed or currently being processed, so that 
\begin{equation}
	f(\vec{x}):=\sum_{i\in D} c_ix_i.
	\label{cost_formulation}
\end{equation}

Let $\theta$ denote a decision-making policy for the MDP. In a general sense, $\theta$ can take actions that depend on the history of the process thus far and might also be randomized (see \cite{Puterman1994}).  The policy is said to be \emph{stationary} and \emph{deterministic} if action $\theta(\vec{x})\in A_\vec{x}$ is always chosen under state $\vec{x}\in S$. The expected long-run average cost \black{per unit time} (or \emph{average cost} for short) under policy $\theta$, given that the initial state of the system is $\vec{x}\in S$, can be expressed as\\
\begin{equation}
	g_\theta(\mathbf{x})=\liminf_{T\rightarrow\infty}T^{-1}\mathbb{E}_\theta \left[\sum_{t=0}^{T-1} f(\mathbf{x}(t))\;\Big |\;\mathbf{x}(0)=\mathbf{x}\right]
	\label{avg_cost}
\end{equation}
where $\vec{x}(t)$ denotes the state of the system at time step $t\in\mathbb{N}_0$. The theory of uniformization implies that under a stationary policy $\theta$, the discretized system has the same average cost $g_\theta(\vec{x})$ as that incurred by its continuous-time counterpart, in which we allow actions to be chosen every time the system transitions from one state to another. \black{In the continuous-time system, after choosing an action $a\in A_\vec{x}$ in state $\vec{x}\in S$, the amount of time until the next transition is exponentially distributed with mean $1/q(\vec{x},a)$, where $q(\vec{x},a)$ is the sum of the transition rates (for example, if the server is processing a job at demand point $i$ then $q(\vec{x},a)=\sum_{j\in D}\lambda_j+\mu_i$). On the other hand, in the uniformized system (with $\Delta=1$), there is a fixed time step of unit length until the next transition, and we also allow `self-transitions' into the same state. It can be seen from (\ref{trans_prob_eqn}) that the probability of a self-transition is $1-q(\vec{x},a)$ and hence the total number of time steps until the next change of state (assuming that the decision-maker continues to choose action $a$) is geometrically distributed with mean $1/q(\vec{x},a)$. Thus, the expected amount of time (and also the cost incurred) until the next change of state is the same in both systems. This provides an intuitive justification for the uniformization method, and we defer to \cite{lippman1975device} and \cite{serfozo1979equivalence} for full technical details.}

The objective of the problem is to minimize the long-run average cost $g_\theta(\vec{x})$; that is, to find a policy $\theta^*$ such that
$$g_{\theta^*}(\vec{x})\leq g_\theta(\vec{x})\;\;\;\;\forall \theta\in \Theta,\;\vec{x}\in S,$$
where $\Theta$ is the set of all admissible policies.

\black{\begin{remark}The MDP model that we describe in this section can be extended without great difficulty to the case of phase-type processing and switching times. In order to consider phase-type processing times, suppose the processing time of a job at node $i$ is distributed as a sum of $k_i\geq 2$ exponentially-distributed service phases, with $\mu_i^{(r)}$ being the service rate for the $r^{\text{th}}$ phase, $r=1,...,k_i$. The state variable $x_i$ should then be the number of service phases remaining to be processed, rather than the number of jobs. Thus, if a new job of type $i$ arrives at a particular time step (with probability $\lambda_i$), this causes $x_i$ to increase by $k_i$. One can define $y_i=\lceil x_i/k_i\rceil$ as the number of jobs remaining (including any that are partially complete). Since jobs are processed in first-come-first-served order, if $x_i=2k_i+1$ (for example) and the server remains at node $i$ at a particular time step then a service phase completion occurs with probability $\mu_i^{(k_i)}$, since there is a partially complete job with only one phase remaining. The cost function should still be based on job counts rather than phase counts, so the single-step cost function is $\sum_{i\in D}c_iy_i$. The uniformization parameter in (\ref{Delta_eqn}) can be adjusted in an obvious way. To incorporate phase-type setup times (between demand points), one can allow the switching rates between adjacent nodes to be edge-dependent and also (if desired) direction-dependent, so that the setup times are distributed as sums of independent exponential random variables with non-identical rates. This is a powerful generalization because phase-type distributions can approximate any continuous distribution to an arbitrary degree of accuracy (\cite{Asmussen2003}).
\label{remark1}\end{remark}}

Given that our MDP has an infinite state space, the average cost $g_\theta(\vec{x})$ can only be finite if the system is stable under $\theta$, in the sense that $\theta$ induces an ergodic Markov chain on $S$.  Let $\rho$ denote the \emph{traffic intensity} of the system, defined as 
$$\rho=\sum_{i\in D} \lambda_i/\mu_i.$$
Our first result establishes that the condition $\rho<1$ is sufficient to ensure the existence of a deterministic stationary policy under which the system is stable.

\begin{theorem}\textbf{(Stability.)} Suppose $\rho < 1$.  Then there exists a deterministic stationary policy $\theta$ such that $g_\theta(\vec{x})=:g_\theta<\infty$ for all $\vec{x}\in S$.\label{stability_theorem}\end{theorem}

\begin{figure}[htbp]
	\begin{center}
    	\includeinkscape[scale = 0.8]{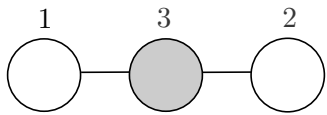}
        \end{center}
	\caption{A network with 2 demand points and 1 intermediate stage. }
	\label{fig2}
\end{figure}

Details of the proof can be found in Appendix \ref{AppA}.  It relies upon results from the literature on polling systems, but these results cannot be applied directly to our system because it is not possible to construct a stationary policy $\theta$ that visits the demand points $i\in D$ in a fixed cyclic pattern and serves each point exhaustively on each visit.  To illustrate this point, consider a network with only two demand points ($D=\{1,2\}$) separated by a single intermediate stage ($N=\{3\}$) as shown in Figure \ref{fig2}. Suppose we wish to implement a simple `polling system' type of policy under which the server moves between the two demand points in an alternating pattern ($1,2,1,2,...$) and, each time it arrives at point $i\in D$, stays there until all jobs have been processed ($x_i=0$) before moving directly to the other point.  Unfortunately, in order to know which action to choose under the state $(3,(0,0))$ (or any other state $\vec{x}$ with $v(\vec{x})=3$), the server must know which demand point was the last to be visited.  Since this information is not included in the system state, the server is forced to follow a nonstationary policy in order to achieve the required alternating pattern.

Our proof circumvents this difficulty by considering a new MDP with a modified state space, in which the system state includes an extra variable indicating the most recent demand point $i$ to have had no jobs present while the server was located there. This enables `polling system' type policies, like the one described above, to be represented as stationary policies. The condition $\rho<1$ is known to ensure that, under an exhaustive `polling'-type policy, the system is stable (\cite{altman1992stability}). We can then use the `approximating sequence' approach developed by Sennott (\cite{sennott1997computation,sennott1999stochastic}) to compute an optimal policy for the modified MDP using value iteration. Finally, we use an inductive argument to show that the decisions taken under the optimal policy for the modified MDP are independent of the extra information that we included in the state space, 
and therefore the optimal policy $\theta^*$ for the modified MDP (under which the system is stable) also qualifies as an deterministic stationary (and optimal) policy for the original MDP. 

It is interesting to note that, by Theorem \ref{stability_theorem} the condition $\rho<1$ is always sufficient to ensure that the system is stable under an exhaustive polling regime, regardless of how long it takes for the server to switch between different demand points. In our model, switching times generally become longer if the switching rate $\tau$ is reduced or extra intermediate stages are inserted in the network.  
Intuitively, if switching times become longer, then the number of jobs present when the server arrives at any particular demand point tends to increase. However, under the exhaustive regime, the server is committed to processing all jobs at a demand point before moving to the next one. Suppose the expected switching times between demand points are large but finite. If the system is unstable, then there must be at least one demand point at which the number of jobs present tends to infinity over time, but this implies that the overall proportion of time that the server spends processing jobs (as opposed to switching between nodes) tends towards 1. This is not possible if $\sum_i \lambda_i/\mu_i<1$, since the proportion of time spent processing jobs at any node $i\in D$ cannot be greater than $\lambda_i/\mu_i$. Thus, regardless of how large the expected switching times are, the system is indeed stable under exhaustive polling when $\rho<1$.

In theory, the `approximating sequences' approach discussed in the proof of Theorem \ref{stability_theorem} can be combined with value iteration in order to compute an optimal policy for our MDP, assuming that $\rho<1$.  In practice, however, the well-known `curse of dimensionality' prevents us from being able to compute optimal policies in systems with more than (roughly) three or four demand points. Therefore we need to develop heuristic methods to obtain easily implementable policies that can achieve strong performances across a range of different possible system configurations. Our proposed heuristic methods are introduced in Section \ref{sec:Index heuristic}.\\

\section{Index heuristics}
\label{sec:Index heuristic}

\text{ }\indent As mentioned in Section \ref{sec:Introduction}, \cite{duenyas1996heuristic} considered a job scheduling problem that has some similarities to ours, although their problem is not defined on a network and does not allow switching or processing times to be interrupted. \black{The approach used in \cite{duenyas1996heuristic} (see also \cite{bell1971rewards,harrison1975priority}) is based on a well-known equivalence between the minimization of congestion-based costs and the maximization of a reward criterion under which the server obtains rewards at a constant rate while it is processing jobs at a particular demand point. We will show that a similar equivalence holds in our problem, despite the inclusion of novel features such as interruptible setup and processing times. In order to motivate the approach, it is necessary to begin by discussing how the problem formulated in Section \ref{sec:Problem Formulation} would be modified if the objective was to minimize expected total \emph{discounted} cost, rather than long-run average cost. Let $\alpha$ be a discount factor satisfying $0<\alpha<1$ (values close to 1 are typically assumed). We will consider a deterministic, stationary policy $\theta$ and use $\theta(\vec{x})$ to denote the action $a\in A_{\vec{x}}$ selected by $\theta$ under state $\vec{x}\in S$. Then the expected total discounted cost under an admissible policy $\theta$ can be expressed as}
\black{\begin{align}V(\vec{x})&=\mathbb{E}_\theta\left[\sum_{t=0}^\infty \alpha^t f(\vec{x}(t))\;\Big |\;\vec{x}(0)=\vec{x}\right]\nonumber\\[8pt]
&=\mathbb{E}_\theta\left[\sum_{t=0}^\infty \alpha^t \sum_{i\in D}c_i[x_i(0)+A_i(t)-B_i(t)]\;\Big |\;\vec{x}(0)=\vec{x}\right],\label{disc_cost}\end{align}}
\black{where $x_i(0)$ is the initial number of jobs at demand point $i\in D$, $A_i(t)$ is the total number of new jobs of type $i$ that arrive by time step $t$ and $B_i(t)$ is the total number of jobs of type $i$ that are completely processed by time step $t$. Since $x_i(0)$ and $A_i(t)$ (for all $i\in D$, $t\geq 0$) are independent of the server's actions, it is clear from (\ref{disc_cost}) that minimizing $V(\vec{x})$ is equivalent to maximizing a weighted sum of expected departure counts, $B_i(t)$. Furthermore, if the server is processing a job of type $i$ at time step $t$, then with probability $\mu_i$ the number of jobs at $i$ is reduced by one and this implies an expected cost saving of $\mu_i\sum_{u=t+1}^\infty \alpha^u c_i=c_i\mu_i\alpha^{t+1}/(1-\alpha)$ compared to a passive policy which allows this job to remain in the system forever. It follows that minimizing $V(\vec{x})$ is equivalent to maximizing}
\black{\begin{equation}W(\vec{x}):=\mathbb{E}_\theta\left[\sum_{t=0}^\infty \alpha^t w(\vec{x}(t),\theta(\vec{x}(t)))\;\Big|\;\vec{x}(0)=\vec{x}\right],\label{disc_reward}\end{equation}}
\black{where the single-step reward function $w(\vec{x},a)$ is defined for $\vec{x}=(v,(x_1,...,x_d))\in S$, $a\in A_{\vec{x}}$ by}
\black{\begin{equation}w(\vec{x},a)=\begin{cases}
\dfrac{\alpha c_i\mu_i}{1-\alpha},&\text{if }v=i\text{ for some }i\in D\text{, }x_i\geq 1\text{ and }a=v,\\[8pt]
0,&\text{otherwise.}\end{cases}\label{disc_reward_function}\end{equation}}
\black{It is also important to note that this equivalence holds under quite general conditions. The fact that setup and processing times are interruptible in our formulation does not cause any complications, as the function $B_i(t)$ in (\ref{disc_cost}) is a step function that increases only when jobs are finally completed and exit from the system.

It is clear from (\ref{disc_reward_function}) that the single-step reward associated with processing a job of type $i$ tends to infinity as $\alpha\rightarrow 1$, because the extra cost associated with holding a job in the system forever becomes infinite. Therefore, following \cite{duenyas1996heuristic}, we derive a bounded `reward rate' by dividing the expected total discounted cost saving from choosing action $a$ under state $\vec{x}$ (namely, $w(\vec{x},a)$) by the discounted length of the  time horizon, $\sum_{t=0}^\infty \alpha^t=(1-\alpha)^{-1}$, to obtain the following reward function for undiscounted problems:}
\black{\begin{equation}r(\vec{x},a)=\lim_{\alpha\rightarrow 1}(1-\alpha)w(\vec{x},a)=\begin{cases}
c_i\mu_i,&\text{if }v=i\text{ for some }i\in D\text{, }x_i\geq 1\text{ and }a=v,\\
0,&\text{otherwise.}\end{cases}\label{reward_formulation}\end{equation}}

\black{\begin{remark}If an extension to phase-type processing times is to be considered (as discussed in Remark \ref{remark1}) one can derive an expression for the reward rate $r(\vec{x},a)$ similar to the one in (\ref{reward_formulation}), except that the condition $x_i\geq 1$ should be replaced by $x_i\equiv 1\; (\emph{mod }k_i)$, since a service completion is possible only if the job at the front of the queue has one phase remaining.\end{remark}}

\black{The reward rates $r(\vec{x},a)$ are useful in the development of index heuristics because they depend only on the server's current location and choice of action. On the other hand, the cost function (\ref{cost_formulation}), despite its simple appearance, has a disadvantage in that the aggregate cost incurred by the system at any given time step must be calculated according to the number of jobs present at every demand point in the system.  Thus, it is more convenient to use reward rates when seeking to develop \emph{index heuristics}, which operate by associating easily computable scores (or \emph{indices}) with the decision options under any given state. Throughout this section, the index heuristics that we develop are based on the reward rates defined in (\ref{reward_formulation}). However, in our numerical experiments in Section \ref{sec:numerical}, we evaluate and compare the performances of different policies according to the original cost formulation presented in Section \ref{sec:Problem Formulation}.}

\subsection{The DVO heuristic}\label{DVO_heuristic}

\text{ }\indent We refer to the heuristic policy proposed in \cite{duenyas1996heuristic} as the `DVO heuristic' for short. The DVO heuristic does not allow processing or switching times to be interrupted, and therefore decision epochs occur only if (1) the server finishes processing a job, (2) the server arrives at a demand point, or (3) the server is idle and a new job arrives in the system. As such, we cannot represent the DVO heuristic as a stationary policy in our system, since the action chosen at a particular time step constrains the actions chosen at future time steps. For example, if the server begins switching from one demand point to another, then it must continue to do so until it arrives at the new demand point. Similarly, if the server remains at a non-empty demand point then it must continue to do so until a job is processed. However, despite its lack of compatibility with our MDP formulation, we can still consider the DVO heuristic as a nonstationary policy in our system and estimate its performance using simulation experiments. We use this approach to provide a useful benchmark for our heuristic policies in Section \ref{sec:numerical}.

The complete steps of the DVO heuristic are described in \cite{duenyas1996heuristic}. For the reader's convenience, we provide a summary of how the heuristic makes decisions at the three different types of decision epoch (1)-(3) mentioned above in Appendix \ref{AppDVO}. 

\subsection{$K$-stop heuristics}\label{sequence_heuristics}

\text{ }\indent Although the DVO heuristic can be applied to our problem and its performance in a particular system can be simulated, it is unlikely to achieve near-optimal performances in general. There are two main reasons for this: (i) it over-constrains the action sets by not allowing switching or processing times to be interrupted; (ii) it chooses the next demand point to switch to without considering the proximity of that demand point to other demand points in the network. The latter of these two properties implies that the heuristic works in a short-sighted (myopic) way, and fails to recognize the potential benefits of moving to demand points that are located near other demand points. This is not a weakness in the problem studied by \cite{duenyas1996heuristic}, since in their model, the switching (or setup) time to move from node $i$ to node $j$ only depends on the destination node $j$, so the currently-occupied node does not have any effect on future switching times. In our network-based formulation, however, it clearly makes sense to consider heuristic policies that can take the network topology into account when making decisions.

In this subsection we introduce a class of heuristic policies designed to exploit the novel features of our problem. We refer to these policies as $K$-stop heuristics, since the decision at any particular time step is made by assuming that the server will visit a sequence of up to $K$ demand points (where $K$ is a pre-determined integer) and serve these demand points until exhaustion. The name `$K$-stop' derives from the fact that the server is assumed to visit a sequence of demand points in the same way that a public transport service visits different `stops' along its route; however, it is important to emphasize that the optimal sequence to be followed by the server is calculated only for the purposes of making a single decision at a particular time step, and the server is not actually committed to following this route in future time steps. Thus, given that each time step is also a decision epoch in our problem formulation, routes are continuously re-optimized and therefore switching and processing times can be interrupted.

Let $K\geq 1$ be a pre-determined integer. At any given time step we consider all non-empty sequences of demand points with length not exceeding $K$; that is, we consider sequences of the form $s=(s_1,s_2,...,s_m)$, where $1\leq m\leq K$ and $s_j\in D$ for each $j=1,2,...,m$. We also require that all elements of the sequence are distinct (that is, $s_j\neq s_k$ for $j,k\in\{1,...,m\}$ with $j\neq k$), and furthermore the first element $s_1$ must be different from the server's current location $v$. For each sequence $s$ we calculate a `reward rate', also referred to as an \emph{index}, and this index is calculated by assuming that the server visits demand points in the order $s_1,...,s_{m}$ and serves each demand point exhaustively before moving to the next one. Some further conditions related to stability and idling (described below) are also used to decide which sequences should be considered `eligible'. After identifying the eligible sequence with the highest index, we choose an action at the current time step by assuming that this sequence should be followed.

Before providing full details of the $K$-stop heuristic, it will be useful to introduce some extra notation. Given a state $\vec{x}=(v,(x_1,...,x_d))$ and a sequence $s$, we define $s_0=v$ for notational convenience; that is, $s_0$ is the server's current node (which may not be a demand point). We also use $\delta(i,j)$ to denote the length of the shortest path (in terms of the number of nodes that must be traversed) from node $i\in V$ to $j\in V$, with $\delta_{ii}:=0$ for $i\in V$. For a given state $\vec{x}$, sequence $s$ and each $j=1,...,|s|$, we define two quantities $T_j(\vec{x},s,t)$ and $R_j(\vec{x},s,t)$ as follows:
\begin{align}&T_j(\vec{x},s,t)=\dfrac{x_{s_j}+\lambda_{s_j}\left[t+\sum_{k=1}^{j-1}\left(\delta(s_{k-1},s_k)/\tau+T_k(\vec{x},s,t)\right)+\delta(s_{j-1},s_j)/\tau\right]}{\mu_{s_j}-\lambda_{s_j}},&t\geq 0,\label{Ts_eqn}\\
&R_j(\vec{x},s,t)=c_{s_j}\mu_{s_j}T_j(\vec{x},s,t),&t\geq 0.\label{Rs_eqn}\end{align}
In words, $T_j(\vec{x},s,t)$ is an approximation for the expected amount of time required for the server to process all jobs at node $s_j$ after arriving there, assuming that it remains idle for $t$ time units at the current node $v$ before beginning to follow the sequence $s$. This approximation is obtained by assuming a fluid-type model of the system dynamics in which jobs arrive at node $i\in D$ at a continuous rate $\lambda_i$ and are processed at a continuous rate $\mu_i$, and the server requires $1/\tau$ time units to move between any adjacent pair of nodes in the network. On the other hand, $R_j(\vec{x},s,t)$ is the total reward earned while the server is processing jobs at node $s_j$ in this fluid model. To make sense of equation (\ref{Ts_eqn}), note that the number of jobs present at node $s_j$ when the server arrives there is given by adding the original number of jobs under state $\vec{x}$ (that is, $x_{s_j}$) to the number of new jobs that arrive while the server is either idling, traveling to one of the earlier demand points in the sequence, processing jobs at one of the earlier demand points, or traveling towards $s_j$.  Furthermore, once the server arrives at node $s_j$, the number of jobs at $s_j$ decreases at a net rate of $\mu_{s_j}-\lambda_{s_j}$. \black{We note that the equations \ref{Ts_eqn} and (\ref{Rs_eqn}) can also be adapted quite easily to more general models with edge-dependent switching rates, implying phase-type setup times, as discussed in Remark \ref{remark1}. In this case, the expression $\delta(s_{k-1},s_k)/\tau$ should be replaced by the minimum possible expected amount of time to travel from node $s_{k-1}$ to $s_k$, taking into account the rates of travel on individual edges.}

It is important to clarify that $t$ is interpreted as a certain amount of `idle time' before the server begins to follow the sequence $s$. During this idle time, the server remains idle at its current node $v$ and (in the event that $v$ is a demand point) does not process any jobs there. This is somewhat contrary to our MDP formulation in Section \ref{sec:Problem Formulation}, which assumes that if the server remains at a non-empty demand point then it must be processing jobs there. However, we make this `idle time' assumption only for the purpose of deriving index quantities for use in our heuristics. Next, for a given state $\vec{x}$ and sequence $s$, we define
\begin{align}&\psi(\vec{x},s,t)=\dfrac{\sum_{k=1}^{|s|} R_k(\vec{x},s,t)}{t+\sum_{k=1}^{|s|}\left[\delta(s_{k-1},s_k)/\tau+T_k(\vec{x},s,t)\right]},\;\;\;\;t\geq 0,\label{psi_eqn}\\[8pt]
&\phi_j(\vec{x},s,t)=\dfrac{\sum_{k=1}^{j} R_{k}(\vec{x},s,t)}{t+\sum_{k=1}^{j}\left[\delta(s_{k-1},s_k)/\tau+T_{k}(\vec{x},s,t)\right]+\delta(s_{j},v)/\tau},\;\;\;\;t\geq 0,\;j\in\{1,...,|s|\}.\label{phi_eqn}\end{align}
We can interpret $\psi(\vec{x},s,t)$ as the average reward per unit time earned while the server follows sequence $s$. On the other hand, $\phi_j(\vec{x},s,t)$ is the average reward earned during a truncated sequence that visits the demand points $s_1,s_2,...,s_j$ and then returns to the starting node $v$. In this case, the average reward calculation includes the time taken to switch back to node $v$. The quantities $\psi(\vec{x},s,t)$ and $\phi_j(\vec{x},s,t)$ are used for slightly different purposes in our heuristics. Recall that we only consider sequences in which the first node $s_1$ differs from the server's current location $s_0$, and therefore the denominators in (\ref{psi_eqn}) and (\ref{phi_eqn}) are always non-zero. 

Before presenting the algorithmic steps for the $K$-stop heuristic, we prove a useful property of the average reward $\psi(\vec{x},s,t)$.
\begin{lemma}For any given state $\vec{x}\in S$ and sequence $s$, the average reward $\psi(\vec{x},s,t)$ is a monotonic function of $t$.
\label{sec3_lem}\end{lemma}

\textit{Proof. }We prove the statement by showing that the derivative $\frac{\partial}{\partial t}\psi(\vec{x},s,t)$ has the same sign (either positive, negative or zero) for all $t\geq 0$. The state variables $x_1,...,x_d$ can be regarded as constants. The system parameters (including arrival rates, processing rates, switching rate, cost rates and the distances between nodes in the network) can also be regarded as constants. Therefore, using a trivial inductive argument, we can show that the quantities $T_j(\vec{x},s,t)$ defined in (\ref{Ts_eqn}) are linear functions of $t$, for $j=1,...,|s|$. Hence, the average rewards $R_j(\vec{x},s,t)$ in (\ref{Rs_eqn}) are also linear in $t$. It follows that the quantity $\psi(\vec{x},s,t)$ in (\ref{psi_eqn}) is a ratio of two linear functions, with the general form $(a_1+b_1t)/(a_2+b_2t)$, where $a_1,a_2,b_1,b_2$ are positive constants. Any such function can be represented graphically as a hyperbola, with a derivative of the form $(a_2b_1-a_1b_2)/(a_2+b_2t)^2$. Hence, the sign of the derivative is the same for any $t\geq 0$. \hfill $\Box$\\

Next, we provide details of the steps used in the $K$-stop heuristic algorithm. Recall that $K\geq 1$ is a pre-determined, fixed integer and we let $\vec{x}=(v,(x_1,...,x_d))$ denote the current state. At each time step, there are three possible cases: (1) the server is at a non-empty demand point, (2) the server is at an empty demand point, (3) the server is at an intermediate stage. The details below explain how actions are chosen in each of these cases.\\ 

\noindent\textbf{$K$-stop heuristic algorithm}
\begin{enumerate}
\item If the server is at a non-empty demand point ($v\in D$ and $x_v>0$) then we perform the following steps:
\begin{enumerate}[(a)]
\item Let $\mathcal S$ be the set of all sequences of the form $s=(s_1,s_2,...,s_m)$, where $1\leq m\leq K$, $s_j\in D$ for each $j\in\{1,2,...,m\}$, $s_1\neq v$ and $s_i\neq s_j$ for any pair of elements $s_i,s_j\in s$ with $i\neq j$. Initialize $\sigma=\emptyset$ \black{as a set of `eligible' sequences from which we will later select a reward-maximizing sequence.}
\item For each sequence $s\in\mathcal S$, define $\beta_j(\vec{x},s,t)$ for $j=1,...,|s|$ as follows:
{\begin{equation}
        \beta_j(\vec{x},s,t):=\begin{cases}\dfrac{\sum_{k=1}^{j} R_k(\vec{x},s,t)}{\sum_{k=1}^{j} T_k(\vec{x},s,t)}\rho+c_v\mu_v(1-\rho),\;\;\;\; &\text{ if }v \notin \{s_1,s_2,...,s_j\},\\[20pt]
	0,\;\;\;\;&\text{ otherwise. }
 \label{beta_eqn}\end{cases}\\[8pt]
    \end{equation}}
If $\frac{\partial}{\partial t}\psi(\vec{x},s,t)\big|_{t=0}\leq 0$ and $\phi_{j}(\vec{x},s,0)\geq\beta_j(\vec{x},s,0)$ for all $j\in\{1,...,|s|\}$, then add $s$ to the set $\sigma$. Otherwise, do not make any change to the set $\sigma$.
\item If $\sigma=\emptyset$, then the action chosen under state $\vec{x}$ should be to remain at node $v$. Otherwise, let $s^*$ denote the sequence in $\sigma$ with the largest value of $\psi(\vec{x},s,0)$, with ties broken according to some fixed priority ordering of the demand points in $D$, \black{applied lexicographically to the sequences in $\sigma$}. The action chosen under $\vec{x}$ should be to switch to the first node on a shortest path from $v$ to $s_1^*$. 
\end{enumerate}
\item If the server is at \black{either an empty demand point ($v\in D$ and $x_v=0$) or an intermediate stage ($v\in N$)} then we perform the following steps:
\begin{enumerate}
\item Let $\mathcal S$ be defined in the same way as in step 1(a). Initialize \black{$\sigma_1=\emptyset$ as a set of `high-priority' eligible sequences and $\sigma_2=\emptyset$ as a set of `low-priority' eligible sequences}.
\item For each sequence $s\in\mathcal S$, add $s$ to $\sigma_2$ if and only if $\frac{\partial}{\partial t}\psi(\vec{x},s,t)\big|_{t=0}\leq 0$. 
\item For each sequence $s\in\sigma_2$, define $\gamma(\vec{x},s,t)$ as follows
\begin{equation}\gamma(\vec{x},s,t):=\dfrac{\sum_{k=1}^{|s|} R_k(\vec{x},s,t)}{\sum_{k=1}^{|s|} T_k(\vec{x},s,t)}\rho,\;\;\;\;t\geq 0.\label{gamma_eqn}\end{equation}
    Let $\vec{y}$ denote a state identical to $\vec{x}$ except that the server is located at $s_1$ (the first node in sequence $s$) instead of $v$. If $\psi(\vec{x},s,0)\geq\gamma(\vec{x},s,0)$ and either (i) $|s|=1$ or (ii) $|s|\geq 2$ and  $\psi(\vec{y},s,0)\geq\gamma(\vec{y},s,0)$, then remove $s$ from $\sigma_2$ and add it to $\sigma_1$. Otherwise, do not make any changes.
\item If $\sigma_1$ is non-empty, set $\sigma=\sigma_1$. Otherwise, set $\sigma=\sigma_2$.
\item Carry out step 1(c).
\end{enumerate}
\end{enumerate}

The steps presented above require some explanation. Consider the first case, where the server is at a non-empty demand point. In this case, $\sigma$ is a set of `eligible' sequences and we choose a sequence from this set that yields the highest average reward, $\psi(\vec{x},s,0)$, in the fluid model. To determine whether some sequence $s\in\mathcal S$ is eligible, we need to check two conditions. If the condition $\frac{\partial}{\partial t}\psi(\vec{x},s,t)\big|_{t=0}\leq 0$ is satisfied, this indicates that we should begin following sequence $s$ immediately (i.e. we should take an immediate step towards the first node in $s$), as the average reward will not increase if we wait for some `idle time' before following $s$. Note that, due to Lemma \ref{sec3_lem}, we know that if $\frac{\partial}{\partial t}\psi(\vec{x},s,t)\big|_{t=0}\leq 0$ then $\frac{\partial}{\partial t}\psi(\vec{x},s,t)\leq 0$ for all $t\geq 0$, so there is no advantage to be gained by waiting for any amount of idle time. The second condition, $\phi_j(\vec{x},s,0)\geq \beta_j(\vec{x},s,0)$, is automatically satisfied if the server returns to its current location $v$ at some point in the first $j$ stops of the sequence (since we define $\beta_j(\vec{x},s,t)=0$ in this case). In other cases, the condition is intended to provide a balance between two important considerations. 
To elaborate on this, note that if $\rho$ is close to 1 then the condition $\phi_j(\vec{x},s,0)\geq \beta_j(\vec{x},s,0)$ is almost equivalent to
$$\dfrac{\sum_{k=1}^{j} T_k(\vec{x},s,0)}{\sum_{k=1}^{j}\left[\delta(s_{k-1},s_k)/\tau+T_k(\vec{x},s,0)\right]+\delta(s_{j},v)/\tau}\geq\rho,$$
which states that if the server serves nodes $s_1,...,s_j$ exhaustively and then returns to node $v$, then the proportion of time spent processing jobs (as opposed to switching between nodes) during this time should be at least $\rho$. This is an important condition for system stability. 
On the other hand, if $\rho$ is close to zero, then $\phi_j(\vec{x},s,0)\geq \beta_j(\vec{x},s,0)$ is almost equivalent to
$$\dfrac{\sum_{k=1}^{j} R_k(\vec{x},s,0)}{\sum_{k=1}^{j}\left[\delta(s_{k-1},s_k)/\tau+T_k(\vec{x},s,0)\right]+\delta(s_{j},v)/\tau}\geq c_v\mu_v,$$
which states that the average reward earned while serving nodes $s_1,...,s_j$ and then returning to $v$ should be at least as great as the average reward that would be earned by continuing to process jobs at the current node; that is, $c_v\mu_v$. Effectively, this states that the server should only move away from node $v$ to process jobs at other demand points if there exists a sequence $s$ for which the average reward is greater than the average reward for remaining at $v$. The condition $\phi_j(\vec{x},s,0)\geq\beta_j(\vec{x},s,0)$ is a generalization of a rule used by the DVO heuristic in \cite{duenyas1996heuristic}, which effectively considers sequences of length 1 only. 

From a computational standpoint, we note that it is not necessary to derive expressions for the derivatives $\frac{\partial}{\partial t}\psi(\vec{x},s,t)$ in terms of the system parameters. Indeed, these expressions become very complicated as the sequence length $|s|$ increases. Instead, we can simply compare the values of $\psi(\vec{x},s,0)$ and $\psi(\vec{x},s,\varepsilon)$, where $\varepsilon$ is some small positive number. This approach is justified by Lemma \ref{sec3_lem}, which implies that $\frac{\partial}{\partial t}\psi(\vec{x},s,t)\big|_{t=0}\leq 0$ if and only if $\psi(\vec{x},s,\varepsilon)\leq \psi(\vec{x},s,0)$.

Next, consider case 2, where the server is at \black{either an empty demand point or an intermediate stage}. In this case, by using two different sets $\sigma_1$ and $\sigma_2$, we effectively separate the eligible sequences into two different subsets, with sequences in $\sigma_1$ being given a higher priority for selection than those in $\sigma_2$. The DVO heuristic in \cite{duenyas1996heuristic} also separates the decision options into two sets when the server is at an empty demand point, but the conditions that we use in our sequence-based algorithm are quite different. Firstly, as in case 1, we introduce a derivative-based condition and require sequence $s$ to satisfy $\frac{\partial}{\partial t}\psi(\vec{x},s,t)\big|_{t=0}\leq 0$ in order to be included in $\sigma_2$. This condition implies that there is no benefit to be gained by waiting for some `idle time' before following $s$. In order to be included in the higher-priority set $\sigma_1$, sequences must also satisfy the condition $\psi(\vec{x},s,0)\geq\gamma(\vec{x},s,0)$ and (for sequences of length greater than one only) $\psi(\vec{y},s,0)\geq\gamma(\vec{y},s,0)$, where $\vec{y}$ is a state identical to $\vec{x}$ except that the server is located at node $s_1$ instead of $v$. We defer discussion of these conditions to the proof of Theorem \ref{pathwise_thm}, where it is shown that they are sufficient to ensure that the sequence $s$ remains included in the set $\sigma_1$ at all stages while the server travels from $v$ to $s_1$, provided that no further jobs arrive in the meantime. This is useful in order to ensure that the server follows a consistent path through the network, rather than changing direction without an obvious reason. We also note that the condition $\psi(\vec{x},s,0)\geq\gamma(\vec{x},s,0)$ is somewhat similar to the condition $\phi_j(\vec{x},s,0)\geq\beta_j(\vec{x},s,0)$ used in case 1, but there are some differences. Firstly, the condition applies to the entire sequence $s$, rather than the subsequences $(s_1,...,s_j)$ for $1\leq j\leq |s|$, so in this sense it is more relaxed than the condition used in case 1. Secondly, the expression for $\gamma(\vec{x},s,t)$ does not include the extra term $c_v\mu_v(1-\rho)$ that can be seen in the expression for $\beta_j(\vec{x},s,t)$, implying that we do not wish to ensure that the average reward obtained by following sequence $s$ is greater than $c_v\mu_v$. Indeed, this is logical since there are no jobs present at node $v$, and therefore it is not possible to earn any immediate reward by remaining there. 

It is worthwhile to emphasize that even if a sequence $s$ fails to satisfy the conditions for inclusion in $\sigma_1$, it may still be included in $\sigma_2$, in which case it may be selected in step 2(e) if $\sigma_1$ is empty. The conditions for inclusion in $\sigma_1$ (that is, $\psi(\vec{x},s,0)\geq\gamma(\vec{x},s,0)$ and $\psi(\vec{y},s,0)\geq\gamma(\vec{y},s,0)$) are primarily intended to promote system stability. However, even if sequence $s$ fails to satisfy these conditions, it may still be better to follow $s$ rather than remaining idle at node $v$. This is intuitive, since idling at an empty demand point \black{or at an intermediate stage} is not necessarily helpful for maintaining stability, so following sequence $s$ should not be seen as a worse option in this regard. 

Our next result states that the $K$-stop heuristic has the property of \emph{pathwise consistency}, which implies that the server follows a consistent path through the intermediate stages of the network as long as the number of jobs in the system remains unchanged. This is an intuitively appealing property, as it means that the server avoids wasting time (by going back and forth between intermediate stages, for example) when moving between demand points.

\begin{theorem}\textbf{(Pathwise consistency.) }Suppose the server is located at an intermediate stage $v\in N$ and the system operates under the $K$-stop heuristic policy. Then there exists a demand point $j^*\in D$ such that the server moves directly along a shortest path to node $j^*$ until either (i) it arrives at node $j^*$, or (ii) a new job arrives in the system.\label{pathwise_thm}\end{theorem}

Proof of Theorem \ref{pathwise_thm} can be found in Appendix \ref{AppB}. It is important to clarify that if a new job arrives in the system while the server is moving towards the demand point $j^*$ referred to in the theorem, then the server may change direction and move towards a different demand point instead. This does not contradict the theorem; indeed, the theorem does not make any claim about what happens after the next job arrival. The theorem essentially states that the server follows a consistent path until the next time a new job arrives, and we are able to use this in order to prove that the expected amount of time until it arrives at a demand point must be finite. We state this as a corollary below and provide a proof in Appendix \ref{AppC}.

\begin{corollary}Suppose the conditions of Theorem \ref{pathwise_thm} apply. Then the expected amount of time until the server arrives at a demand point is finite. More specifically, if $T_{\textup{switch}}$ denotes the amount of time until the server arrives at a demand point, then
$$\mathbb{E}[T_{\textup{switch}}]\leq\frac{M}{\tau}\left(\frac{\Lambda+\tau}{\tau}\right)^{2(M-1)},$$
where $\Lambda:=\sum_{i=1}^d \lambda_i$ and $M:=\max_{\{i\in N,\;j\in D\}}\delta(i,j)$ denotes the maximum distance between an intermediate stage and a demand point.\label{pathwise_corollary}\end{corollary}

\black{The next result concerns a special case of the problem in which job types are homogeneous, which means that the arrival rates $\lambda_i$, service rates $\mu_i$ and holding costs $c_i$ are identical for all $i\in D$. In this scenario, we are able to prove that the system is stable under the policy given by the $K$-stop heuristic. Furthermore, if $V$ is a complete graph (so that all demand points are directly connected to each other) then the policy given by the $K$-stop heuristic is optimal. We emphasize that these statements hold for any number of demand points and any $K\geq 1$, and the stability part of the result also holds for any network layout.} 

\black{\begin{theorem}Suppose we have a homogeneous system in which $\lambda_1=...=\lambda_d$, $\mu_1=...=\mu_d$ and $c_1=...=c_d$.  Then, for any $K\geq 1$, the system is stable under the $K$-stop heuristic policy when $\rho<1$. Furthermore, if $V$ is a complete graph then the $K$-stop heuristic policy is optimal.\label{homogeneous_thm}\end{theorem}}

\black{Proof of the theorem can be found in Appendix \ref{App_homsystem}. Subject to the conditions of the theorem, it can be shown that the server visits all demand points infinitely often under the $K$-stop policy and also serves demand points exhaustively on each visit. This ensures that the system is stable. In the case where $V$ is a complete graph, it can also be shown that the $K$-stop policy directs the server to switch to the demand point with the largest number of jobs if it is currently at an empty demand point, and the same rule is used by an optimal policy. In the latter scenario, the optimal policy is a type of `Serve the Longest Queue' (SLQ) policy, and we note that SLQ policies have previously been studied in certain kinds of polling systems. In particular, \cite{liu1992polling} showed that a SLQ policy is optimal for a `symmetric' polling system, similar to the one described in Theorem \ref{homogeneous_thm}, in which arrival and service rates are the same at all demand points. However, to the best of our knowledge, no similar result has previously been proved for a system with interruptible switching and processing times.}

It should be noted that if $K\geq 2$ then the number of sequences in $\mathcal S$ increases \black{rapidly} with the number of demand points $d$, implying that the computational requirements of the $K$-stop heuristic become unmanageable in large-scale problems. 
In the next subsection we propose an alternative heuristic under which the number of indices to be calculated at any time step increases only linearly with $d$, enabling greater scalability.

\subsection{($K$ from $L$)-stop heuristic}\label{rpoint}

\text{ }\indent As explained in Section \ref{sequence_heuristics}, the $K$-stop heuristic considers all possible sequences of demand points $(s_1,s_2,...,s_m)$ at each time step, for each $1\leq m\leq K$. The only restrictions are that all demand points in the sequence are distinct and the first demand point must be different from the server's current location. This implies that the number of sequences to be considered at each time step is of order $d^K$, which grows polynomially with the number of demand points $d$ (for fixed $K$) and grows exponentially with $K$. Hence, in order for the $K$-stop heuristic to be computationally feasible, $d$ and $K$ must be relatively small. In our numerical experiments in Section \ref{sec:numerical} we restrict attention to systems with $d\leq 8$ and consider $K\in\{1,2,3,4\}$.

In this subsection we propose a modified version of the $K$-stop heuristic in which the number of sequences considered at each time step increases only linearly with $d$. Note that if $K=1$ then the $K$-stop heuristic already has this property, since it considers only sequences of the form $(j)$, for $j\in D$. Our modified heuristic works as follows: suppose the system is in state $\vec{x}\in S$ at an arbitrary time step. First, we carry out the same steps as if we are using the $K$-stop heuristic with $K=1$, and calculate the indices $\psi(\vec{x},(j),0)$ for each $j\in D$. In this step we also allow the demand point $j$ to be equal to the server's current location $v$ (unlike in the standard $K$-stop heuristic) and set $\psi(\vec{x},(v),0)$ equal to $c_v\mu_v$ if $v$ is non-empty, and zero otherwise. After computing $\psi(\vec{x},(j),0)$ for each $j\in D$, we then form a set $\mathcal L$ of size $L$ (where $L\leq d$ is a pre-determined integer) consisting of a limited number of demand points. We consider a couple of different ways of selecting the demand points to be included in $\mathcal L$:\\

\begin{itemize}
\item \textbf{Impartial method: }In this method, we simply choose the $L$ demand points with the highest indices $\psi(\vec{x},(j),0)$, regardless of their positions in the network. (Ties are broken arbitrarily.)
\item \textbf{Stratified method: }In some systems, there may be a natural way of dividing the network into `clusters' of demand points, with any pair of demand points in the same cluster being relatively close to each other. In this case, we can select a pre-determined number of demand points from each cluster (taking the ones with the highest indices), in such a way that the total number of demand points selected is $L$.\\
\end{itemize}

In the impartial method we also enforce the following rule: if the server is at an empty demand point or an intermediate stage then we divide the demand points into two sets. The first set consists of sequences $(j)$ such that $\psi(\vec{x},(j),0)\geq \gamma(\vec{x},(j),0)$ and the second set consists of sequences $(j)$ such that $\psi(\vec{x},(j),0)<\gamma(\vec{x},(j),0)$, where $\gamma(\cdot)$ is defined in (\ref{gamma_eqn}). If the first set includes at least $L$ sequences then we select the $L$ sequences with the highest values of $\psi(\vec{x},(j),0)$ to be included in $\mathcal L$. Otherwise, all of the sequences in the first set are included in $\mathcal L$ and we obtain the remaining sequences by choosing the sequences with the highest indices in the second set. This is consistent with the prioritization rule described in step 2(c) of the $K$-stop heuristic. In the stratified method a similar rule is used, except the two sets are formed for each cluster separately, and in each cluster we choose a pre-determined number of demand points from the two sets, with the first set being given priority over the second set as in the impartial method.

After forming the set $\mathcal L$, we then consider all possible sequences $(s_1,...,s_m)$ for $1\leq m\leq K$, where $s_j\in \mathcal L$ for each $j=1,...,m$, and carry out the rest of the steps in the $K$-stop heuristic as described in Section \ref{sequence_heuristics}. The sequences are required to satisfy the same eligibility conditions described in Section \ref{sequence_heuristics} in order to be selected.

We note that the impartial method is simpler and might often perform better than the stratified method, but the stratified method offers a potential advantage in that it forces the server to consider moving to other clusters in the network. Under the impartial method, there is a risk that if all of the $L$ demand points selected are in close proximity to the server's current location then the server acts in a short-sighted way, as (given that demand points are selected based on indices for sequences of length one only) it fails to detect the potential benefits of moving to another cluster and serving multiple demand points within that cluster.

\begin{figure}[htbp]
        \begin{center}
	   \vspace{0.5cm} 
	   \includeinkscape[scale = 0.8]{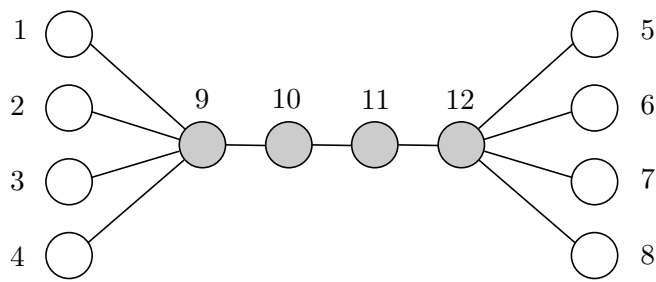}
        \end{center}
	\caption{A network with 4 demand points on the left, 4 demand points on the right and 4 intermediate stages. }
	\label{sec33_fig}
\end{figure}

As an example, consider the system shown in Figure \ref{sec33_fig}, with 8 demand points and 4 intermediate stages. Suppose we use the ($K$ from $L$)-stop heuristic with $K=2$ and $L=4$. For each demand point $j=1,...,8$ we calculate the index $\psi(\vec{x},(j),0)$. Under the impartial method, assuming that the server is at a non-empty demand point, we choose the 4 demand points $j$ with the highest indices and then consider all possible sequences of length $1$ or $2$ involving these 4 demand points only. The total number of sequences to consider is $4+(4!/2!)=16$. If the server is at an empty demand point or an intermediate stage then the process is similar except we prioritize sequences $(j)$ that satisfy the condition $\psi(\vec{x},(j),0)\geq \gamma(\vec{x},(j),0)$. In this case we still obtain $16$ sequences. Under the stratified method, a logical approach (given the layout of the network) is to define $\{1,2,3,4\}$ as one cluster and $\{5,6,7,8\}$ as another cluster, and select two demand points from each according to their index values. Suppose, for example, we choose 2 and 3 from the first cluster and 6 and 7 from the second cluster. Then, in the next step, we consider all possible sequences of length 1 or 2 consisting of demand points from the set $\{2,3,6,7\}$. Thus, we again consider 16 sequences in total.

In general, in large systems (with a lot of demand points), the bulk of the computational effort required is in the calculation of indices $\psi(\vec{x},(j),0)$ for each $j\in D$. Following this, the number of sequences to be considered is $\sum_{m=1}^K (L!/(L-m)!)$, which is independent of $d$ and can be kept relatively small. Thus, if $K$ and $L$ are fixed then the computational effort required by the ($K$ from $L$)-stop heuristic increases only linearly with $d$. If $L=d$ then the ($K$ from $L$)-stop heuristic becomes equivalent to the $K$-stop heuristic. Thus, it is obvious that the ($K$ from $L$)-stop heuristic should perform worse than the $K$-stop heuristic in general, since it considers a smaller number of possible sequences. However, in the next section, we show that it may be able to achieve a similar performance at smaller computational expense.\\

\section{Numerical results}\label{sec:numerical}

\text{ }\indent In this section we report the results of numerical experiments in order to compare the performances of the heuristics described in Section \ref{sec:Index heuristic}. \black{In Section \ref{sec41} we focus on a specific network layout, with two `clusters' of demand points separated by a series of intermediate stages. In Section \ref{sec42} we present results from problem instances with randomly-generated network layouts. Finally, in Section \ref{sec43} we provide additional details of the computational requirements of our heuristics.}

\subsection{Two clusters of demand points}\label{sec41}

\begin{figure}[htbp]
        \begin{center}
	   \vspace{0.5cm}
	   \includeinkscape[scale = 0.8]{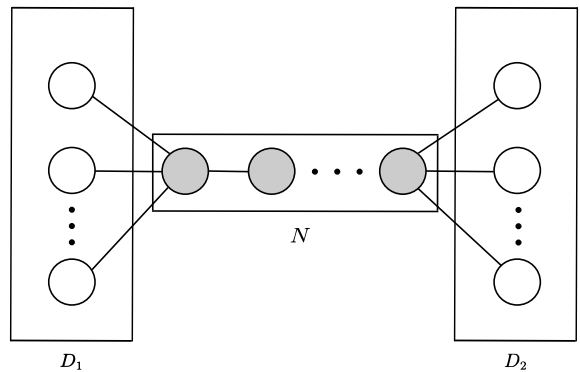}
        \end{center}
	\caption{A diagrammatic representation of the network with $d_1$ demand points on the left, $d_2$ demand points on the right and $n$ intermediate stages.}
	\label{fig4}
\end{figure}

In this subsection we consider a relatively simple network layout, shown in Figure \ref{fig4}, in which two distinct `clusters' of demand points are separated by a chain of $n$ intermediate stages. The demand points on the left-hand side belong to a cluster denoted by $D_1$ of size $d_1$, and similarly the demand points on the right-hand side form a cluster $D_2$ of size $d_2$. In order to move from $D_1$ to $D_2$ or vice versa, the server must pass through all of the intermediate stages in succession. Also, as the figure indicates, in order to move from one demand point to another point in the same cluster, it must pass through one intermediate stage (it is not possible to move directly from one demand point to another). By adjusting the value of $n$ we can vary the distance between the two clusters. In the case $n=1$, we obtain a special case where all demand points are equidistant from each other. In these experiments we consider $1\leq n\leq 6$, $1\leq d_1\leq 4$ and $1\leq d_2\leq 4$. Thus, the number of demand points $d$ satisfies $2\leq d\leq 8$.

Our numerical study is based on $10,000$ randomly-generated problem instances. In each instance we uniformly sample the values of $d_1$, $d_2$ and $n$ from the ranges specified above. We also randomly generate the values of $\rho$, $\tau$, $\lambda_i$, $\mu_i$ and $c_i$ (for each $i\in D$) using a method that ensures consideration of a wide range of different scenarios for the system parameters. 
For full details of our parameter generation methods, please refer to Appendix \ref{AppD}.

In each problem instance, we test the performances of the heuristic policies in Section \ref{sec:Index heuristic} using simulation experiments. The DVO heuristic described in Section \ref{DVO_heuristic} can be regarded as a nonstationary policy for our MDP and its performance can be simulated. This provides a useful benchmark for our other heuristics. We test the $K$-stop heuristic for each $K\in\{1,2,3,4\}$. We note that in instances where $d<K$ (i.e. the number of demand points is smaller than $K$), the $K$-stop heuristic becomes equivalent to the $(K-1)$-stop heuristic, as it is not possible to define a sequence of length $K$ such that all demand points in the sequence are distinct (but the definition of the heuristic allows it to consider sequences of length smaller than $K$). For example, if $d=2$ then the $2$-stop, $3$-stop and $4$-stop heuristics are equivalent. We also implement the ($K$ from $L$)-stop heuristic with $K=2$ and $L=4$ in each instance, using both the impartial and stratified methods. In the stratified case, we adopt the obvious strategy of defining $D_1$ and $D_2$ as separate clusters. For implementation details of our simulation experiments, 
please refer to Appendix \ref{AppE}.

In instances where $d$ is small, it may be possible to compute the optimal long-run average cost $g^*$ using dynamic programming (specifically, relative value iteration). Although DP algorithms require a finite state space, the `approximating sequence' method of \cite{sennott1997computation} can be used to obtain the infinite-state optimal value as a limit of finite-state optimal values (for further explanation, see the proof of Theorem \ref{stability_theorem}). In each problem instance we have used an iterative method, described in Appendix \ref{AppF}, to test whether or not it is computationally feasible to obtain $g^*$ using DP. We have found that it is usually only feasible to compute $g^*$ if $d\leq 3$ and $\rho$ is not too large. 
In total, we have been able to compute $g^*$ in 1229 of the 10,000 instances. 

Table \ref{overall performance for RVIA} shows, for each heuristic policy, a 95\% confidence interval for the mean percentage suboptimality of the heuristic in comparison to the optimal value $g^*$ computed using relative value iteration, based on the results from 1229 `small' problem instances. For each heuristic policy $\theta\in\{\text{DVO, 1-stop, 2-stop, 3-stop}\}$, the percentage suboptimality is calculated as $100\times(g_\theta-g^*)/g^*$. The 10th, 25th, 50th, 75th and 90th percentiles of the distribution of percentage suboptimality for each heuristic are also reported. Given that $d\leq 3$ in all of these small instances, the 4-stop policy is equivalent to the 3-stop policy, so we do not include it in the table. Similarly, the (2 from 4)-stop heuristics are equivalent to the 2-stop heuristic because they always consider all demand points in the network as potential sequence elements, so we do not include them either. The table shows that the DVO heuristic (which does not allow switching or processing times to be interrupted) is about 22\% suboptimal on average in these instances. The $K$-stop heuristics are able to improve upon the DVO heuristic very significantly. The $K=2$ and $K=3$ policies are within 4\% of optimality on average, and in more than 50\% of instances they are within 3\%. As expected, the performance tends to improve as $K$ increases (at the expense of greater computation time), although increasing $K$ from 2 to 3 gives only a small additional improvement.

\begin{table}[htbp] 
	\centering
	\caption{\black{Percentage suboptimalities of heuristic policies in 1229 `small' problem instances.}}
	\label{overall performance for RVIA}
	\resizebox{0.7\textwidth}{!}{
		\begin{tabular}{lcccccccccccccccccccc} 
			\toprule 
			Heuristic &Mean&10th pct.&25th pct.&50th pct.&75th pct.&90th pct.\\
			\midrule
			DVO & $22.05\pm0.83$ & 6.61 & 11.68 & 19.50 & 28.88 & 40.03 \\
                1-stop & $5.97\pm0.45$ & 0.73 & 1.92 & 3.93 & 7.20 & 13.08 \\
                2-stop & $3.96\pm0.24$ & 0.42 & 1.49 & 2.92 & 5.01 & 8.44 \\
                3-stop & $3.73\pm0.24$ & 0.35 & 1.39 & 2.75 & 4.64 & 7.80 \\
			\bottomrule 
	    \end{tabular}}
\end{table}

Table \ref{overall performance for DVO} shows a comparison between the heuristic policies across all $10,000$ instances, with the DVO heuristic used as a benchmark. In each instance we report 95\% confidence intervals for the mean percentage improvements of the $K$-stop and ($K$ from $L$)-stop heuristics over the DVO heuristic, and also report the 10th, 25th, 50th, 75th and 90th percentiles of the distributions for the percentage improvement. We find that all of the $K$-stop and ($K$ from $L$)-stop heuristics that we consider are able to improve significantly over the DVO heuristic, with the mean improvements ranging from about 9\% to 11\%. Clearly, the ability to interrupt switching and processing times offers major advantages, as the decision-maker is able to respond more rapidly to the arrivals of new jobs. We also observe again that the performance of the $K$-stop heuristic tends to improve as $K$ increases, although for $K\geq 3$ the marginal extra improvements become quite small. While larger values of $K$ enable more long-sighted choices of actions, it is not necessarily clear that these should result in dramatic improvements over smaller $K$ values, as the server is always able to change its direction at each time step and is never committed to following a sequence through to its end. Indeed, the larger $K$ is, the less likely it becomes that the server follows a sequence through to completion. Nevertheless, the results do seem to indicate a trend for larger $K$ values to yield stronger policies. It is also encouraging to note that both versions of the (2 from 4)-stop heuristic yield performances very close to that of the 2-stop heuristic. Recall that, in general, the ($K$ from $L$)-stop heuristic is supposed to be a more computationally scalable version of the $K$-stop heuristic. We accept a small loss of performance in exchange for greater scalability. The results indicate that the simpler impartial version of the (2 from 4)-stop heuristic tends to perform slightly better than the stratified version, indicating that there is no obvious benefit in forcing the server to consider moving to a different cluster in the network. 

\begin{table}[htbp] 
	\centering
	\caption{\black{Percentage improvements of heuristic policies with respect to the DVO policy in 10,000 problem instances.}}
	\label{overall performance for DVO}
	\resizebox{0.8\textwidth}{!}{
		\begin{tabular}{lcccccccccccccccccccc} 
			\toprule 
			Heuristic &Mean&10th pct.&25th pct.&50th pct.&75th pct.&90th pct.\\
			\midrule
		1-stop & $9.48\pm0.16$ &0.17 & 4.69 & 9.45 & 14.59 & 19.35\\
            2-stop & $10.50\pm0.15$ & 1.69 & 5.64 & 10.24 & 15.31 & 20.08\\
            (2 from 4)-stop [imp.] & $10.48\pm0.15$ & 1.68 & 5.63 & 10.24 & 15.31 & 20.05\\
            (2 from 4)-stop [str.] & $10.30\pm0.15$ & 1.42 & 5.43 & 10.08 & 15.24 & 19.97\\
            3-stop & $10.92\pm0.14$ & 2.25 & 6.05 & 10.63 & 15.73& 20.23\\
            4-stop & $11.14\pm0.14$ & 2.76 & 6.41 & 10.85 & 15.79 & 20.29\\
			\bottomrule 
	    \end{tabular}}
\end{table}

\black{Next, we investigate the effects of the system parameters on these results by categorizing the $10,000$ problem instances according to the values of specific parameters. Three parameters of particular interest to us are the number of intermediate stages $n$, the traffic intensity $\rho=\sum_i \lambda_i/\mu_i$ and the relative switching rate, which we define as $\eta:=\tau/(\sum_i \lambda_i)$ (so that it indicates the relative speed of switching compared to the frequency of new job arrivals). Table \ref{improvements with varying n} summarizes the relative performances of the heuristics for each $n\in\{1,2,3,4,5,6\}$. Note that the first row of the table shows the improvements of the 1-stop policy over the DVO heuristic, and the remaining rows show the additional improvements of the other heuristics over the 1-stop heuristic. We have chosen to present the results in this way in order to neatly summarize the benefits of allowing switching and processing times to be interrupted (shown in the first row of the table) and also the additional benefits of allowing the heuristic policies to make longer-sighted decisions (shown in the remaining rows). Table \ref{traffic intensity rho} shows similar results for $\rho$, and Table \ref{various tau} shows similar results for $\eta$. Note that our method for randomly generating the system parameter values (detailed in Appendix \ref{AppD}) implies that, in any problem instance, $\rho$ is equally likely to fall within any of the 4 intervals shown as columns in Table \ref{traffic intensity rho}, and similarly $\eta$ is equally like to be fall within any of the 6 intervals shown as columns in Table \ref{various tau}.}

\begin{table}[htbp]
	\centering
	\caption{Percentage improvements of the 1-stop policy (vs. the DVO heuristic) and the 2-stop, (2 from 4)-stop, 3-stop and 4-stop policies (vs. the 1-stop policy) for different values of $n$.}
	\label{improvements with varying n}
	\resizebox{\textwidth}{!}{ 
		\begin{tabular}{lccccccccccccccccccccccccc} 
		\toprule 
		& $n = 1$ & $n = 2$ & $n = 3$ & $n = 4$ & $n = 5$ & $n = 6$\\
            Heuristic & [1671 instances] & [1642] & [1654] & [1706] & [1663] & [1664]\\
		\midrule
            1-stop (vs. DVO)&  $10.58\pm0.50$ & $10.71\pm0.36$ & $10.18\pm0.36$ & $9.14\pm0.36$ & $8.46\pm0.38$ & $7.86\pm0.39$\\
            2-stop (vs. 1-stop)& $1.56\pm0.15$ & $0.88\pm0.10$ & $0.80\pm0.13$ & $0.76\pm0.14$ & $1.04\pm0.16$ & $1.15\pm0.19$\\
            (2 from 4)-stop [imp.] (vs. 1-stop)&  $1.45\pm0.14$ & $0.85\pm0.10$ & $0.86\pm0.12$ & $0.82\pm0.13$ & $1.05\pm0.16$ & $1.08\pm0.19$\\
            (2 from 4)-stop [str.] (vs. 1-stop)&  $1.31\pm0.13$ & $0.78\pm0.09$ & $0.65\pm0.12$ & $0.56\pm0.17$ & $0.77\pm0.19$ & $0.88\pm0.20$\\
			3-stop (vs. 1-stop)& $2.09\pm0.17$ & $1.28\pm0.11$ & $1.28\pm0.13$ & $1.25\pm0.15$ & $1.38\pm0.18$ & $1.52\pm0.20$\\
			4-stop (vs. 1-stop)& $2.23\pm0.19$ & $1.51\pm0.12$ & $1.46\pm0.14$ & $1.44\pm0.15$ & $1.63\pm0.18$ & $1.86\pm0.20$\\
			\bottomrule 
	\end{tabular}}
\end{table}

The first row of Table \ref{improvements with varying n} appears to suggest a trend for the improvements given by the 1-stop policy (compared to the DVO heuristic) to decrease as $n$ increases. To make sense of this, it is useful to bear in mind that under any policy and any network design, the simulated proportion of time that the server spends processing jobs should be approximately equal to $\rho$, assuming that the system is stable. Thus, regardless of how small or large $n$ is, the server should spend approximately the same proportion of time visiting the intermediate stages of the network. \black{In the $n=1$ case, the server derives the maximum possible advantage from being able to switch direction, as it can switch from the intermediate stage to any of the demand points in a single transition.} Under the DVO heuristic, on the other hand, the server has no ability to change direction. It is also worthwhile to note that the `impartial' version of the (2 from 4)-stop heuristic appears to outperform the `stratified' version across all values of $n$. Intuition might suggest that the stratified version should become stronger as $n$ increases, as the clusters in the network become more distinct in this situation and the stratified version is designed to ensure that demand points in both clusters are always considered as potential destinations for switching. However, this intuition is not borne out by the results. It appears that allowing the demand points in the set $\mathcal L$ to be selected solely according to the indices given by the 1-stop policy (regardless of their locations in the network) is consistently the most effective approach. 

\begin{table}[htbp]
	\centering
	\caption{Percentage improvements of the 1-stop policy (vs. the DVO heuristic) and the 2-stop, (2 from 4)-stop, 3-stop and 4-stop policies (vs. the 1-stop policy) for different values of $\rho$.}
	\label{traffic intensity rho}
	\resizebox{0.8\linewidth}{!}{ 
		\begin{tabular}{lccccccccccccccccccccccccc} 
			\toprule 
			 &$0.1\leq \rho < 0.3$ & $0.3\leq \rho < 0.5$ & $0.5\leq \rho < 0.7$ & $0.7\leq \rho < 0.9$\\
              Heuristic & [2448 instances] & [2514] & [2498] & [2540]\\
			\midrule
            1-stop (vs. DVO)& $16.64\pm0.25$ & $11.26\pm0.25$ & $7.00\pm0.27$ & $3.27\pm0.26$\\
            2-stop (vs. 1-stop)& $0.55\pm0.09$ & $0.79\pm0.11$ & $1.37\pm0.13$ & $1.40\pm0.13$\\
            (2 from 4)-stop [imp.] (vs. 1-stop)& $0.57\pm0.09$ & $0.83\pm0.11$ & $1.33\pm0.13$ & $1.33\pm0.13$\\
            (2 from 4)-stop [str.] (vs. 1-stop)& $0.53\pm0.09$ & $0.70\pm0.11$ & $1.17\pm0.13$ & $0.90\pm0.17$\\
			3-stop (vs. 1-stop)& $0.87\pm0.09$ & $1.26\pm0.12$ & $1.84\pm0.15$ & $1.88\pm0.15$\\
			4-stop (vs. 1-stop)& $1.05\pm0.10$ & $1.45\pm0.13$ & $2.09\pm0.15$ & $2.14\pm0.15$\\
			
			\bottomrule 
	\end{tabular}}
\end{table}

\black{Next, we discuss the effect of the traffic intensity, $\rho$.} By examining the first row in Table \ref{traffic intensity rho}, we observe that the improvements given by the 1-stop policy (vs. the DVO heuristic) are much greater when $\rho$ is small than when it is large. Indeed, if $\rho$ is small, then the server spends a relatively large proportion of its time traversing the intermediate stages of the network. In such situations, the ability to react immediately to the arrival of a new job (by changing the direction of travel) clearly offers major advantages. On the other hand, when $\rho$ is large, queue lengths tend to become longer at the demand points and the server spends more of its time processing jobs at the demand points. It is useful to bear in mind that the DVO heuristic retains the ability to make new decisions every time the server finishes processing a job (which happens often when $\rho$ is large), whereas it does not have the ability to make new decisions while the server is switching, so in this sense larger $\rho$ values work in its favor. Therefore, it is more difficult for the 1-stop policy to gain an advantage when $\rho$ is large. The remaining rows in Table \ref{traffic intensity rho} show that, in most cases, the improvements of the $K$-stop and ($K$ from $L$)-stop policies (for $K\geq 2$) over the 1-stop policy tend to increase as $\rho$ increases. This makes sense intuitively, as the longer-sighted heuristics are able to assess the congestion levels at multiple demand points when choosing a sequence to follow. In heavy traffic systems, it is more important to recognize situations where switching to a different part of the network allows the server to process jobs at a number of demand points in succession. 

\begin{table}[htbp]
	\centering
	\caption{Percentage improvements of the 1-stop policy (vs. the DVO heuristic) and the 2-stop, (2 from 4)-stop, 3-stop and 4-stop policies (vs. the 1-stop policy) for different values of $\eta$.}
	\label{various tau}
	\resizebox{\textwidth}{!}{ 
		\begin{tabular}{lccccccccccccccccccccccccc} 
			\toprule 
		 &$0.1\leq \eta < 0.4$ & $0.4\leq \eta < 0.7$ & $0.7\leq \eta < 1$ & $1\leq \eta < 4$ & $4\leq \eta < 7$ & $7\leq \eta < 10$\\
            Heuristic & [1628 instances] & [1688] & [1599] & [1737] & [1680] & [1668]\\
			\midrule 
            1-stop (vs. DVO) & $7.52\pm0.47$ & $9.01\pm0.42$ & $9.67\pm0.39$ & $11.58\pm0.38$ & $9.96\pm0.36$ & $9.03\pm0.33$\\
            2-stop (vs. 1-stop) & $1.56\pm0.20$ & $1.43\pm0.17$ & $1.39\pm0.16$ & $0.99\pm0.12$ & $0.53\pm0.10$ & $0.32\pm0.08$\\
            (2 from 4)-stop [imp.] (vs. 1-stop)& $1.53\pm0.19$ & $1.44\pm0.16$ & $1.41\pm0.15$ & $0.99\pm0.12$ & $0.48\pm0.10$ & $0.27\pm0.10$\\
            (2 from 4)-stop [str.] (vs. 1-stop)& $1.21\pm0.22$ & $1.17\pm0.17$ & $1.12\pm0.17$ & $0.77\pm0.15$ & $0.44\pm0.10$ & $0.27\pm0.08$\\
			3-stop (vs. 1-stop)& $2.36\pm0.21$ & $2.02\pm0.19$ & $1.94\pm0.17$ & $1.28\pm0.13$ & $0.76\pm0.10$ & $0.49\pm0.09$\\
			4-stop (vs. 1-stop)& $2.83\pm0.22$ & $2.38\pm0.19$ & $2.22\pm0.18$ & $1.43\pm0.13$ & $0.86\pm0.10$ & $0.47\pm0.10$\\
			\bottomrule 
	\end{tabular}}
\end{table}

\black{Finally, we discuss the effect of the switching rate parameter, $\eta$.} 
From Table \ref{various tau} we observe that the mean percentage improvement of the 1-stop policy (vs. the DVO heuristic) tends to decrease as $\eta$ becomes small, but also as $\eta$ becomes large. Indeed, when $\eta$ is small, switching between nodes is relatively slow and the 1-stop policy is more likely to be deterred from changing direction, as this will result in too much wasted time traversing the intermediate stages of the network. On the other hand, when $\eta$ is large, switching between nodes is relatively fast and in this situation it becomes less likely that any arrivals occur while the server is traversing the intermediate stages, so the server has no reason to change its course. 
\black{In general, the improvement over the DVO policy should be greater when direction changes are frequent.} From the remaining rows in Table \ref{various tau}, we observe that the $K$-stop and ($K$ from $L$)-stop heuristics (for $K\geq 2$) are able to achieve greater improvements over the 1-stop policy when $\eta$ is small. Indeed, a primary motivation for using the longer-sighted heuristics is that they can recognize the effects of distances between the demand points and plan a sequence of visits accordingly. When $\eta$ is large, distances become less important as the server is always able to switch quickly between any two demand points and react quickly to new job arrivals, so it becomes more difficult for the longer-sighted heuristics to achieve improvements over the 1-stop policy. 

\black{\subsection{Randomly-generated network layouts}\label{sec42}}

\begin{figure}[htbp]
    \centering
    \begin{tikzpicture}[scale=0.8]
        \foreach \x in {-5,-4,...,5} {
          \ifnum\x=0
          \else
            \draw[gray, very thin] (\x,-5) -- (\x,-1);
            \draw[gray, very thin] (\x,1) -- (\x,5);
          \fi
        }
        
        \foreach \y in {-5,-4,...,5} {
          \ifnum\y=0
          \else
            \draw[gray, very thin] (-5,\y) -- (-1,\y);
            \draw[gray, very thin] (1,\y) -- (5,\y);
          \fi
        }

        \draw[thick,-] (0,-5.5) -- (0,5.5); 
        \draw[thick,-] (-5.5,0) -- (5.5,0) node[right] {}; 

        \foreach \x in {-5,-4,...,5} {
            \ifnum\x=0\else \draw (\x,0) -- (\x,-0); \fi
        }
        \foreach \y in {-5,-4,...,5} {
            \ifnum\y=0\else \draw (0,\y) -- (-0,\y); \fi
        }

        \begin{scope}[shift={(0,0)}]
            \draw[fill=white, draw=black] (1,2) circle (8pt);  
            \draw[fill=white, draw=black] (3,4) circle (8pt);  
            \draw[fill=white, draw=black] (5,3) circle (8pt);  
    
            \foreach \x in {1,2,3,4} {
                \draw[fill=lightgray, draw=black] (\x,3) circle (8pt);
            }
    
            \draw (1,2.28) -- (1,2.72);             
            \draw (1.28,3) -- (1.72,3);             
            \draw (2.28,3) -- (2.72,3);             
            \draw (3.28,3) -- (3.72,3);             
            \draw (4.28,3) -- (4.72,3);             
            \draw (3,3.28) -- (3,3.72);             

        \end{scope}

        \begin{scope}[shift={(-5,0)}]
            \draw[fill=white, draw=black] (1,5) circle (8pt);  
            \draw[fill=white, draw=black] (1,4) circle (8pt);  
            \draw[fill=white, draw=black] (1,1) circle (8pt);  
            \draw[fill=white, draw=black] (2,2) circle (8pt);  
            \draw[fill=white, draw=black] (3,5) circle (8pt);  
            \draw[fill=white, draw=black] (4,4) circle (8pt); 
    
            \draw[fill=lightgray, draw=black] (1,2) circle (8pt);  
            \draw[fill=lightgray, draw=black] (1,3) circle (8pt);  
            \draw[fill=lightgray, draw=black] (2,3) circle (8pt);  
            \draw[fill=lightgray, draw=black] (2,4) circle (8pt);  
            \draw[fill=lightgray, draw=black] (2,5) circle (8pt);  
            \draw[fill=lightgray, draw=black] (3,4) circle (8pt); 
    
            \draw (1,1.28) -- (1,1.72);             
            \draw (1,2.28) -- (1,2.72);             
            \draw (1,3.28) -- (1,3.72);             
            \draw (1,4.28) -- (1,4.72);             
            \draw (2,2.28) -- (2,2.72);             
            \draw (2,3.28) -- (2,3.72);             
            \draw (2,4.28) -- (2,4.72);             
            \draw (3,4.28) -- (3,4.72);             

            \draw (1.28,2) -- (1.72,2);             
            \draw (1.28,3) -- (1.72,3);             
            \draw (1.28,4) -- (1.72,4);             
            \draw (1.28,5) -- (1.72,5);             
            \draw (2.28,4) -- (2.72,4);             
            \draw (2.28,5) -- (2.72,5);             
            \draw (3.28,4) -- (3.72,4);             
        \end{scope}

        \begin{scope}[shift={(-5,-5)}]
            \draw[fill=white, draw=black] (0,2) circle (8pt);  
            \draw[fill=white, draw=black] (0,3) circle (8pt);  
            \draw[fill=white, draw=black] (2,4) circle (8pt);  
            \draw[fill=white, draw=black] (4,0) circle (8pt);  
            \draw[fill=white, draw=black] (4,3) circle (8pt);  

            \draw[fill=lightgray, draw=black] (1,2) circle (8pt);  
            \draw[fill=lightgray, draw=black] (1,3) circle (8pt);  
            \draw[fill=lightgray, draw=black] (2,2) circle (8pt);  
            \draw[fill=lightgray, draw=black] (2,3) circle (8pt);  
            \draw[fill=lightgray, draw=black] (3,2) circle (8pt); 
            \draw[fill=lightgray, draw=black] (3,3) circle (8pt);  
            \draw[fill=lightgray, draw=black] (4,1) circle (8pt);  
            \draw[fill=lightgray, draw=black] (4,2) circle (8pt); 
    
            \draw (0,2.28) -- (0,2.72);             
            \draw (1,2.28) -- (1,2.72);             
            \draw (2,2.28) -- (2,2.72);             
            \draw (3,2.28) -- (3,2.72);             
            \draw (4,2.28) -- (4,2.72);             
            \draw (2,3.28) -- (2,3.72);             
            \draw (4,1.28) -- (4,1.72);             
            \draw (4,0.28) -- (4,0.72);             

            \draw (0.28,2) -- (0.72,2);             
            \draw (1.28,2) -- (1.72,2);             
            \draw (2.28,2) -- (2.72,2);             
            \draw (3.28,2) -- (3.72,2);             
            
            \draw (0.28,3) -- (0.72,3);             
            \draw (1.28,3) -- (1.72,3);             
            \draw (2.28,3) -- (2.72,3);             
            \draw (3.28,3) -- (3.72,3);             
        \end{scope}

        \begin{scope}[shift={(0,-5)}]
            \draw[fill=white, draw=black] (1,1) circle (8pt); 
            \draw[fill=white, draw=black] (2,0) circle (8pt);  
            \draw[fill=white, draw=black] (2,3) circle (8pt);  
            \draw[fill=white, draw=black] (2,4) circle (8pt); 
            \draw[fill=white, draw=black] (3,1) circle (8pt);  
            \draw[fill=white, draw=black] (3,3) circle (8pt); 
            \draw[fill=white, draw=black] (5,0) circle (8pt);  
            \draw[fill=white, draw=black] (5,3) circle (8pt); 
    
            \draw[fill=lightgray, draw=black] (2,1) circle (8pt);  
            \draw[fill=lightgray, draw=black] (2,2) circle (8pt);  
            \draw[fill=lightgray, draw=black] (3,0) circle (8pt);  
            \draw[fill=lightgray, draw=black] (3,2) circle (8pt);  
            \draw[fill=lightgray, draw=black] (4,0) circle (8pt);  
            \draw[fill=lightgray, draw=black] (4,1) circle (8pt);  
            \draw[fill=lightgray, draw=black] (4,2) circle (8pt);  
            \draw[fill=lightgray, draw=black] (4,3) circle (8pt);  
            \draw[fill=lightgray, draw=black] (5,1) circle (8pt); 
            \draw[fill=lightgray, draw=black] (5,2) circle (8pt); 
            
            \draw (2,3.28) -- (2,3.72);             
            
            \draw (2,2.28) -- (2,2.72);             
            \draw (3,2.28) -- (3,2.72);             
            \draw (4,2.28) -- (4,2.72);             
            \draw (5,2.28) -- (5,2.72);             
            
            \draw (2,1.28) -- (2,1.72);             
            \draw (3,1.28) -- (3,1.72);             
            \draw (4,1.28) -- (4,1.72);             
            \draw (5,1.28) -- (5,1.72);             
            
            \draw (2,0.28) -- (2,0.72);             
            \draw (3,0.28) -- (3,0.72);             
            \draw (4,0.28) -- (4,0.72);             
            \draw (5,0.28) -- (5,0.72);             

            \draw (2.28,0) -- (2.72,0);             
            \draw (3.28,0) -- (3.72,0);             
            \draw (4.28,0) -- (4.72,0);             
            
            \draw (1.28,1) -- (1.72,1);             
            \draw (2.28,1) -- (2.72,1);             
            \draw (3.28,1) -- (3.72,1);             
            \draw (4.28,1) -- (4.72,1);             

            \draw (2.28,2) -- (2.72,2);             
            \draw (3.28,2) -- (3.72,2);             
            \draw (4.28,2) -- (4.72,2);             

            \draw (2.28,3) -- (2.72,3);             
            \draw (3.28,3) -- (3.72,3);             
            \draw (4.28,3) -- (4.72,3);             
            
        \end{scope}

    \end{tikzpicture}
    \caption{\black{4 randomly-generated network layouts with demand points shown in white and intermediate stages shown in gray, after removal of redundant intermediate stages.}}
    \label{fig:four_different_graphs}
    \label{fig_randnetworks}
\end{figure}
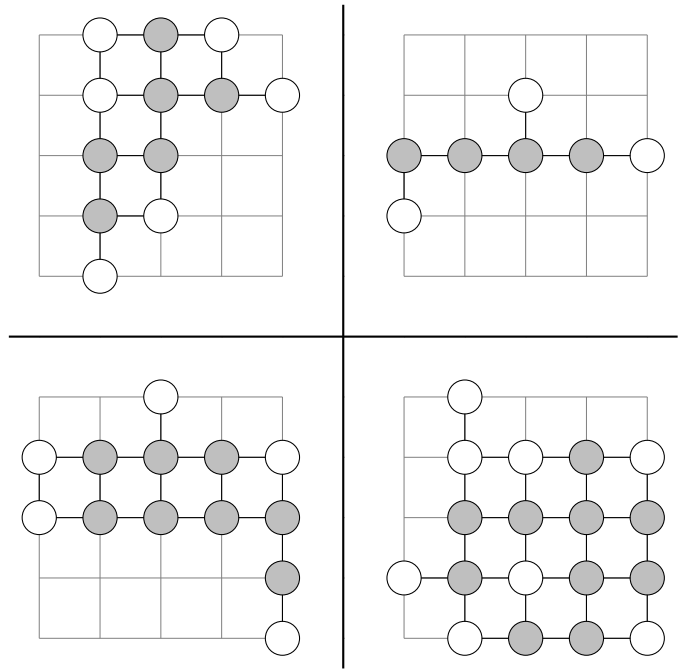

\black{The results that we reported in Section \ref{sec41} were based on a specific network layout, depicted in Figure \ref{fig4}. In order to investigate how well our heuristics perform on other network topologies, we have also carried out $3000$ additional experiments based on randomly-generated network layouts. In these experiments we used the same method for generating the values of arrival rates, processing rates and switching rates as in Section \ref{sec41} (described in Appendix \ref{AppD}) but changed the method for generating the network layout. To generate a random network layout, we begin by creating a $5\times 5$ integer lattice in which $25$ nodes are connected via horizontal and vertical edges. We then randomly select $d$ of these nodes, where $2\leq d\leq 8$, and designate these as demand points. Hence, the position of any demand point $i\in D$ can be represented by a pair of coordinates $(a_i,b_i)$ with $a_i,b_i\in\{1,2,3,4,5\}$, and the shortest path between two demand points $i,j\in D$ is the Manhattan distance $|a_i-a_j|+|b_i-b_j|$.  The remaining nodes in the lattice are designated as intermediate stages, although some of these stages may be redundant (in the sense that they will never be visited under a sensible policy), and therefore can be eliminated from the network. Figure \ref{fig_randnetworks} shows 4 of the random network layouts generated in our experiments.}

\black{As in Section \ref{sec41}, we have used relative value iteration to compute exact suboptimalities of our heuristics in `small' instances where this is feasible. In total, 634 of the 3000 instances satisfied the feasibility criteria described in Appendix \ref{AppF}, and these suboptimality percentages are shown in Table \ref{rand_networks_smallinstances}. By comparing this table with Table \ref{overall performance for RVIA} from earlier, we may observe that all of the heuristics have slightly higher mean suboptimality percentages than in the previously-considered network design. From a qualitative point of view, however, the results are similar, with the $K$-stop heuristic tending to perform better as $K$ increases. It is also interesting to note that the $50^{\text{th}}$ percentiles for the heuristics are similar to those in Table \ref{overall performance for RVIA}, and even lower in some cases. This suggests that the increases in the mean suboptimalities are due to worse tail performances; in other words, there may be some randomly-generated networks on which our heuristics perform particularly poorly. One possible explanation for this is that, in the randomly-generated networks, there may sometimes be multiple shortest paths between pairs of demand points, and currently our heuristics select a shortest path based on an arbitrary priority ordering of the nodes, without taking into account of which path(s) might leave the server in a better position if a switch between demand points is interrupted. Selection of the `best' shortest path (which might depend on the current system state) would significantly complicate our heuristics, and could be a direction for future work.}

\begin{table}[htbp] 
	\centering
	\color{black}
	\caption{\black{Percentage suboptimalities of heuristic policies in 634 `small' problem instances on randomly-generated networks.}}
	\label{randomly overall performance for RVIA}
	\resizebox{0.7\textwidth}{!}{
		\begin{tabular}{lcccccccccccccccccccc} 
			\toprule 
			Heuristic &Mean&10th pct.&25th pct.&50th pct.&75th pct.&90th pct.\\
			\midrule
			DVO &$26.56\pm1.95$ & 5.82 & 10.66 & 21.11 & 33.85 & 50.89\\
                1-stop & $7.46\pm1.47$ & 0.82 & 1.77 & 3.46 & 7.35 & 15.02\\
                2-stop & $4.73\pm0.53$ & 0.62 & 1.43 & 2.86 & 5.01 & 9.49\\
                3-stop & $4.53\pm0.51$ & 0.54 & 1.36 & 2.76 & 4.95 & 9.19\\
			\bottomrule 
	    \end{tabular}}
\label{rand_networks_smallinstances}\end{table}

\black{Table \ref{rand_networks_allinstances} shows the percentage improvements achieved by the $K$-stop and ($K$ from $L$)-stop heuristics against the DVO policy in all of the 3000 instances with randomly-generated networks. The table can be interpreted in the same way as Table \ref{overall performance for DVO} from Section \ref{sec41}, except that we have excluded the `stratified' version of the ($2$ from $4$)-stop heuristic, because it was already shown to perform worse than the impartial version of the heuristic in the previous results, and additionally the division of demand points into `clusters' does not work in such an obvious way in the random networks. By comparing Table \ref{rand_networks_allinstances} with Table \ref{overall performance for DVO} from earlier, we may observe that the improvements of our heuristics compared to the DVO policy are slightly greater than in the previously-considered network design. As discussed earlier, the main weaknesses of the DVO policy are that it works in a short-sighted way and does not allow interruptions of setup or processing times. Thus, the results suggest that the benefits of interruptions are greater in the randomly-generated networks. Indeed, it can be seen from the examples in Figure \ref{fig_randnetworks} that, while moving from one demand point to another, the server might pass near (or through) another demand point, in which case it might be advantageous to switch to that point.}

\begin{table}[htbp] 
	\centering
	\color{black}
	\caption{\black{Percentage improvements of heuristic policies with respect to the DVO policy in 3000 problem instances on randomly-generated networks.}}
	\label{randomly overall performance for DVO}
	\resizebox{0.8\textwidth}{!}{
		\begin{tabular}{lcccccccccccccccccccc} 
			\toprule 
			Heuristic &Mean&10th pct.&25th pct.&50th pct.&75th pct.&90th pct.\\
			\midrule
		1-stop & $12.75\pm0.36$ & 3.42 & 7.81 & 12.70 & 18.09 & 22.95\\
            2-stop & $13.95\pm0.36$ & 4.40 & 8.99 & 14.20 & 19.44 & 24.17\\
            (2 from 4)-stop [imp.] & $13.94\pm0.39$ & 4.43 & 8.92 & 14.15 & 19.49 & 24.24\\
            3-stop & $14.34\pm0.36$ & 4.72 & 9.27 & 14.60 & 19.75 & 24.65\\
            4-stop & $14.45\pm0.36$ & 4.79 & 9.35 & 14.68 & 19.85 & 24.71\\
			\bottomrule 
	    \end{tabular}}
\label{rand_networks_allinstances}\end{table}

\black{We have also investigated the effects of individual parameters (such as $\rho$ and $\eta$) on the results for the randomly-generated networks and found that, from a qualitative point of view, the effects are similar to those discussed in Section \ref{sec41}, so we have omitted these for brevity.}

\black{\subsection{Computational requirements of our heuristics}\label{sec43}}

\black{An advantage of using index-based heuristics, as proposed in our paper, is that they have much lighter computational requirements than (for example) a reinforcement learning algorithm, or any method that relies on extensive training in order to learn a strong decision-making policy. However, as discussed previously, the computational requirements of the $K$-stop heuristic increase rapidly as $K$ increases (or as the number of demand points increases), due to the growth in the size of $\mathcal S$, the set of candidate sequences. This provides motivation for using the ($K$ from $L$)-stop heuristics that we propose in Section \ref{rpoint}. }

\black{In Table \ref{Av_run_times} we present a summary of the average computation times (in seconds) needed by the various heuristics to select an action at a single time step. Our experiments were performed on 108 problem instances with randomly-generated networks and 8 demand points ($d=8$). The software used was Python 3.7.13, with the PyPy just-in-time compiler (\url{http://pypy.org}) used to speed up computations, and all experiments were carried out on an Apple M1 Pro (8-core CPU) with 16 GB unified memory, running macOS Sequoia. The times shown in the table are average times (in seconds) needed to carry out all decision-making steps at a single discrete time step. For example, in the case of the $K$-stop heuristic, this includes the time needed to construct the set of sequences $\mathcal S$ and then select an action according to the steps described in the algorithm as presented in Section \ref{sequence_heuristics}. The running times are averaged over all time steps within each instance and then averaged over all instances. The results show that all of the heuristics are able to make decisions within 0.001 seconds (and much less in some cases), which implies that they are suitable for use in fast-changing systems with hundreds or even thousands of state changes per second. As expected, the increase in running time as $K$ increases is significant, with the average running time tending to increase by a factor of between 4 and 10 each time $K$ increases by one. Notably, the ($2$ from $4$)-stop heuristic has a very short running time and is even faster than the $2$-stop heuristic because, once the subset of 4 demand points has been selected (which requires only a short amount of time), the number of sequences of length $2$ to be evaluated is much smaller than in the $2$-stop heuristic.}

\begin{table}[htbp] 
	\centering
	\color{black}
	\caption{\black{Average running times (in seconds per time step) for action selection by the DVO policy, 1-stop, 2-stop, (2 from 4)-stop [imp.], 3-stop, and 4-stop heuristics, over 108 problem instances with 8 demand points on randomly-generated networks.}}
	\label{Av_run_times}
	\resizebox{0.9\textwidth}{!}{
		\begin{tabular}{lcccccccccccccccccccc} 
			\toprule 
			Heuristic &DVO &1-stop &2-stop &(2 from 4)-stop [imp.] &3-stop &4-stop\\
			\midrule
            Av.run.times & $7.19 \times 10^{-5}$ & $1.72 \times 10^{-6}$ & $7.21 \times 10^{-6}$ & $3.17 \times 10^{-6}$ & $6.36 \times 10^{-5}$ & $4.61 \times 10^{-4}$ \\
			\bottomrule 
	    \end{tabular}}
\end{table}

\section{Conclusions}
\label{sec:conclusions}

\text{ }\indent The main novelty of the job scheduling problem studied in this paper lies in its network-based formulation, which allows setup and processing times to be interruptible and also enables the modeling of complex dependence structures between the setup requirements of different tasks. We consider a highly stochastic, infinite-horizon problem in which jobs of different types arrive at random points in time and the server's setup and processing times are also random. The dynamic nature of our problem implies that decision epochs occur very frequently, and this creates some challenges. For example, as discussed in Section \ref{sec:Problem Formulation}, we are unable to directly leverage results from the literature on polling systems in order to prove the existence of a deterministic, stationary policy under which the system is stable, because polling-type policies are nonstationary under our MDP formulation. However, the stability result can be established (provided that $\rho<1$) by proving the equivalence of a similar MDP in which the system state includes extra information.

The index policies that we develop in Section \ref{sec:Index heuristic} are influenced by the heuristic approaches used in \cite{duenyas1996heuristic}, but these approaches must be adapted in order to exploit the novel features of our problem. In particular, our network-based formulation motivates the use of a long-sighted, sequence-based algorithm that takes the topology of the network into account when making decisions. In addition, we ensure that the algorithm makes new decisions at each time step, so that setup and processing times can always be interrupted. We also introduce derivative-based conditions in the steps of the $K$-stop algorithm and use these to show that the resulting policies possess the property of `pathwise consistency', which ensures that the server always proceeds to a demand point in finite time. Furthermore, we propose a modified version of the algorithm (known as the ($K$ from $L$)-stop algorithm) that scales much more readily to systems with many demand points. \black{In special cases of the problem we can prove system stability and optimality of our heuristics (Theorem \ref{homogeneous_thm}), but in more general cases we must study their performance empirically.} 

The numerical results in Section \ref{sec:numerical} demonstrate the advantages of using heuristics that are well-tailored to the novel features of our problem. In small problem instances, we are able to show that the $K$-stop and ($K$ from $L$)-stop heuristics are much closer to optimality than the unmodified DVO heuristic. In larger systems, we are also able to observe the benefits of allowing the server to make long-sighted decisions and to change its course of action when necessary. 
From a practical perspective, it is encouraging to see that the impartial version of the (2 from 4)-stop heuristic performs almost as well as the more computationally intensive 2-stop heuristic, which in turn offers considerable improvements over the more myopic 1-stop heuristic. Thus, in larger problem instances, the (2 from 4)-stop heuristic may be seen as a strong candidate to achieve significant cost savings over simpler alternatives, without incurring excessive computational costs.\\

\noindent\textbf{Acknowledgements}

The first author was funded by a PhD studentship from the Engineering and Physical Sciences Research Council (EPSRC), under grant EP/W523811/1.\\

\newpage
\Large \textbf{Online Appendices}\\
\normalsize

\begin{appendices}

\section{Proof of Theorem \ref{stability_theorem}.}\label{AppA}

\textit{Proof. }Consider a modified MDP in which the state space is 
\begin{equation}\tilde{S}:=\left\{(v,w,(x_1,...,x_d))\;|\;v\in V,\;w\in D,\;x_i\geq 0 \text{ for }i\in D\right\},\label{modified_state_space}\end{equation}
where the extra variable $w$ represents the most recent demand point at which the server witnessed no jobs present. More precisely, if $w=i$ then this indicates that at a certain time step $t_0$ in the history of the process the system was in a state with $v=i$ and $x_i=0$, and none of the states visited in the more recent time steps (between $t_0$ and the present time) had the server at another demand point $j\in D\setminus \{i\}$ at which there were no jobs present.  (We can set $w$ to an arbitrary value when the process is initialized.)  All other aspects of the modified MDP formulation (e.g. actions and costs) are the same as in the original version.  The transitions of the process do not lose their memoryless property when $w$ is included, since the knowledge that $w=i$ at a particular time step is sufficient to specify the probability distribution for its value at the next step; specifically, $w$ is guaranteed to remain unchanged unless either of the following two cases applies:
\begin{enumerate}
\item The server is at a demand point $j\in D\setminus \{i\}$ with $x_j=1$ and chooses action $j$, in which case there is a probability of $\mu_j$ that we have $w=j$ at the next time step and a probability of $1-\mu_j$ that we still have $w=i$.  
\item The server is at an intermediate stage $k\in N$ adjacent to a demand point $j\in D\setminus \{i\}$ with $x_j=0$ and chooses action $j$, in which case there is a probability of $\tau$ that we have $w=j$ at the next time step and a probability of $1-\tau$ that we still have $w=i$.  
\end{enumerate}
Consider a `polling system' policy $\theta^{[P]}$ under which the server visits the demand points in a repeating sequence $(1,2,...,d, 1, 2, ..., d, 1, 2, ...)$ and, upon arriving at any demand point $i$, remains there until the number of jobs has been reduced to zero ($x_i=0$) before moving to the next demand point in the sequence.  If we already have $x_i=0$ when the server arrives at node $i$, then the server immediately moves to the next demand point.  It should be noted that, since switches require an exponentially distributed amount of time, it may happen that the server decides to move away from demand point $i$ when $x_i=0$ but a new job arrives at $i$ while the server is still located at $i$.  In this case, under our proposed policy, the server continues trying to move to the next demand point rather than processing the new job.  We assume that the server always chooses the shortest path (in terms of the number of intermediate stages that must be traversed) between two demand points and, in the case where two or more paths are tied for the shortest length, the path is selected according to some fixed priority ordering of the nodes. This policy can be represented as a stationary policy in our modified MDP by specifying actions according to the simple rule that if $w=i\leq d-1$ then the server attempts to move towards demand point $i+1$ by taking the next step on the shortest path to that node; or, if it is already at $i+1$, then it remains there.  Similarly if $w=d$ then the server attempts to move to demand point $1$ (or remains there).  

Under the proposed policy $\theta^{[P]}$, the system behaves as a polling system with an exhaustive polling regime (meaning that demand points are served until they are empty).  Given that $\rho<1$, Lemma 3.1 in \cite{altman1992stability} implies that the system is stable and there exists a probability distribution $\{\pi_{\theta^{[P]}}(\vec{x})\}_{\vec{x}\in \tilde{S}}$ such that $\pi_{\theta^{[P]}}(\vec{x})$ is the long-run proportion of time spent in state $\vec{x}\in \tilde{S}$.  In the next part of the proof we show that it is possible to use value iteration to compute an optimal policy for the modified MDP.  Given that the state space $S$ is infinite, this requires the use of the `approximating sequences' method developed in \cite{sennott1999stochastic}. Let $\Psi$ denote the modified MDP with state space $\tilde{S}$ defined in (\ref{modified_state_space}) and let $(\Psi_0,\Psi_1,\Psi_2,...)$ denote a sequence of MDPs that are defined on finite spaces, so that the MDP $\Psi_m$ (for $m\in\mathbb{N}_0$) has state space $\tilde{S}_m$ given by
$$\tilde{S}_m:=\left\{(v,w,(x_1,...,x_d))\;|\;v\in V,\;w\in D,\;0\leq x_i\leq m \text{ for }i\in D\right\}.$$
Thus, in the MDP $\Psi_m$, we do not allow the number of jobs at any node $i\in D$ to be greater than $m$.  We do this by modifying the transition probabilities so that if $\vec{x}$ is a state with $x_i=m$ for some $i\in D$, then the arrival of a new job at node $i$ is impossible and instead the `self-transition' probability $p_{\vec{x},\vec{x}}(a)$ is increased by $\lambda_i$.  Let $\theta_m^*$ be an optimal policy for the MDP $\Psi_m$.  We can show that the sequence $(\theta_0^*,\theta_1^*,...)$ converges to an optimal policy for the infinite-state MDP $\Psi$, but this requires certain conditions to be verified.  Specifically, we must show that the assumptions (AC1)-(AC4) described on p. 169 of \cite{sennott1999stochastic} hold for the sequence $(\Psi_0,\Psi_1,\Psi_2,...)$. (See also \cite{sennott1997computation}, pp. 117-118 for an equivalent set of assumptions.)

Assumption (AC1) in \cite{sennott1999stochastic} states that, for the finite-state MDP $\Psi_m$, there exists a constant $g_m$ and a function $h_m:\tilde{S}_m\rightarrow\mathbb{R}$ satisfying
\begin{equation}g_m+h_m(\vec{x})=f(\vec{x})+\min_{a\in A_\vec{x}}\left\{\sum_{\vec{y}\in\tilde{S}}p_{\vec{x},\vec{y}}(a)h_m(\vec{y})\right\}\;\;\;\;\;\;\forall \vec{x}\in \tilde{S}_m.\label{opt_eqs_n}\end{equation}
We can easily show that, given any two states $\vec{x},\vec{y}\in\tilde{S}_m$, there exists a stationary policy $\theta$ such that $\vec{y}$ is accessible from $\vec{x}$ in the Markov chain induced by $\theta$.  This can be achieved by considering arbitrary states $\vec{x},\vec{y}\in\tilde{S}_m$ and identifying a sequence of random transitions that would cause the system to transition from $\vec{x}$ to $\vec{y}$ under policy $\theta$ (which can be specified differently for each pair of states $(\vec{x},\vec{y})$).  It then follows from \cite{Puterman1994} (p. 478) that $\Psi_m$ belongs to the `communicating' class of multichain MDP models, and the constant $g_m$ in (\ref{opt_eqs_n}) is the optimal long-run average cost for $\Psi_m$.  Moreover, the values $h_m(\vec{x})$ for states $\vec{x}\in\tilde{S}_m$ can be computed using the well-known method of value iteration, in which we define $v_m^{(k)}(\vec{x})$ as the optimal (minimal) total cost in a finite-horizon problem with $k$ stages initialized in state $\vec{x}$ and then compute $h_m(\vec{x})=\lim_{k\rightarrow\infty}(v_m^{(k)}(\vec{x})-v_m^{(k)}(\vec{z}))$, with the reference state $\vec{z}\in\tilde{S}_m$ chosen arbitrarily. 

In order to verify assumptions (AC2)-(AC4) we will need to use several properties of the functions $v_m^{(k)}(\vec{x})$.  Let $\vec{x}^{i+}$ denote a state identical to $\vec{x}$ except that one extra job is present at demand point $i\in D$. The required properties are:
\begin{align}&\text{1. }v_m^{(k)}(\vec{x})\leq v_m^{(k+1)}(\vec{x})\;\;\;\;\forall\;m,k\in\mathbb{N}_0,\;\vec{x}\in\tilde{S}_m.\label{stab_prop1}\\
&\text{2. }v_m^{(k)}(\vec{x})\leq v_m^{(k)}(\vec{x}^{i+})\;\;\;\;\forall\;m,k\in\mathbb{N}_0,\;i\in D,\;\vec{x}\in\tilde{S}_m\text{ such that }x_i<m.\label{stab_prop2}\\
&\text{3. }v_m^{(k)}(\vec{x})\leq v_{m+1}^{(k)}(\vec{x})\;\;\;\;\forall\;m,k\in\mathbb{N}_0,\;\vec{x}\in\tilde{S}_m.\label{stab_prop3}\\
&\text{4. Fix }m\in\mathbb{N}_0\text{ and let }\vec{x}=(v,w,(x_1,...,x_d))\text{ and }\vec{x}'=(v',w',(x_1',...,x_d'))\text{ be two states in }\tilde{S}_m\nonumber\\
&\;\;\;\text{ with }v=v'\text{ and }x_i=x_i'\text{ for }i=1,...,d,\text{ but }w\neq w'.\text{ Then }v_m^{(k)}(\vec{x})=v_m^{(k)}(\vec{x}')\text{ for }k\in\mathbb{N}_0.\label{stab_prop4}\end{align}
All of the properties above are logical and can be proved using induction on $k$. We have omitted details of the induction arguments in order to avoid making this proof excessively long, but they are quite straightforward and only require some care in considering the different possible actions that might be chosen by the relevant finite-horizon optimal policies. Property (\ref{stab_prop1}) states that the optimal expected total cost is increasing with the number of stages remaining, $k$. Property (\ref{stab_prop2}) states that this cost is increasing with the number of jobs initially present at any demand point. Property (\ref{stab_prop3}) states that this cost is increasing with the maximum number of jobs, $m$, allowed to be present at any demand point. Finally, property (\ref{stab_prop4}) states that the variable $w$ has no effect on the optimal expected total cost, which makes sense as it does not impose any constraints on the actions that may be chosen in the $k$ remaining stages.

We proceed to verify (AC2)-(AC4). Assumption (AC2) states that $\lim\sup_{m\rightarrow\infty} h_m(\vec{x})<\infty$ for each $\vec{x}\in\tilde{S}$.  Consider the polling-type policy $\theta^{[P]}$ described earlier. It is clear that the Markov chain induced by $\theta^{[P]}$ has a unichain structure on the state space $\tilde{S}$, since the state $\vec{z}=(1,1,(0,0,...,0))$ is accessible from any other state under this policy (this can be seen from the fact that, in between consecutive visits to demand point $1$, it is always possible for no new jobs to arrive at any demand points). Let $J_{\theta^{[P]}}(\vec{x},\vec{z})$ denote the expected total cost incurred until the system enters state $\vec{z}$, given that it is initialized in state $\vec{x}$ and follows policy $\theta^{[P]}$. Since $\vec{z}$ is positive recurrent, it follows from standard theory (see \cite{sennott1999stochastic}, pp. 298-302) that $J_{\theta^{[P]}}(\vec{x},\vec{z})<\infty$ for all $\vec{x}\in\tilde{S}$. Next, for $m\in\mathbb{N}_0$, let $J_{m,\theta^{[P]}}(\vec{x},\vec{z})$ be defined in an analogous way to $J_{\theta^{[P]}}(\vec{x},\vec{z})$ except that we consider applying the policy $\theta^{[P]}$ to the finite-state MDP $\Psi_m$ instead of the infinite-state MDP $\Psi$. It can easily be shown that $J_{m,\theta^{[P]}}(\vec{x},\vec{z})<J_{\theta^{[P]}}(\vec{x},\vec{z})$ for all $m\in\mathbb{N}$, since the amount of time that the server spends at any demand node (and, hence, the cost incurred) is stochastically smaller under $\Psi_m$ than under $\Psi$. We now follow similar arguments to those in the proof of Proposition 8.2.1, p. 171 in \cite{sennott1999stochastic}. By using property (\ref{stab_prop1}) and the fact that $v_m^{(k)}(\vec{x})$ is defined as the expected $k$-stage cost under an optimal policy, we have 
\begin{equation}v_m^{(k)}(\vec{x})\leq v_m^{(p)}(\vec{x})\leq v_{m,\theta}^{(p)}(\vec{x})\;\;\;\;\forall \vec{x}\in\tilde{S},\;m\in\mathbb{N}_0,\;p\geq k,\label{AC2_ineq2}\end{equation}
where $v_{m,\theta}^{(k)}$ is the expected $k$-stage cost under an arbitrary policy $\theta$. Let $\theta$ be defined to mimic policy $\theta^{[P]}$ until the system reaches state $\vec{z}$ and then follow an optimal finite-horizon policy for $k$ steps. Then, from (\ref{AC2_ineq2}) it follows that $v_m^{(k)}(\vec{x})\leq J_{m,\theta^{[P]}}(\vec{x},\vec{z})+v_m^{(k)}(\vec{z})$. Hence, using the previous arguments, we have
$$v_m^{(k)}(\vec{x})-v_m^{(k)}(\vec{z})\leq J_{m,\theta^{[P]}}(\vec{x},\vec{z})\leq J_{\theta^{[P]}}(\vec{x},\vec{z})\;\;\;\;\forall m\in\mathbb{N}_0,\;\vec{x}\in\tilde{S}.$$
Since $h_m(\vec{x})=\lim_{k\rightarrow\infty}(v_m^{(k)}(\vec{x})-v_m^{(k)}(\vec{z}))$, this establishes that the sequence $\{h_m(\vec{x})\}_{m\in\mathbb{N}}$ is bounded above and hence $\lim\sup_{m\rightarrow\infty} h_m(\vec{x})<\infty$ as required.

Assumption (AC3) states that there exists a constant $Q\geq 0$ such that $-Q\leq \lim\inf_{m\rightarrow\infty} h_m(\vec{x})$ for all $\vec{x}\in\tilde{S}$. 
By using property (\ref{stab_prop2}) and taking limits as $k\rightarrow\infty$, we obtain
\begin{equation}h_m(\vec{x})\leq h_m(\vec{x}^{i+})\;\;\;\; \forall i\in D,\;\vec{x}\in\tilde{S}_m\text{ such that }x_i<m.\label{AC3_ineq2}\end{equation}
This shows that the function $h_m$ attains a minimum on the subset of states with no jobs present, which we denote as $U$. That is,
$$\arg\min_{\vec{x}\in\tilde{S}_
m}h_m(\vec{x})\in \left\{(v,w,(x_1,...,x_d))\;|\;v\in V,\;w\in D,\;x_i=0 \text{ for all }i\in D\right\}=:U.$$
Let $\vec{u}^*$ be a state that attains the minimum above. We will assume that $\vec{u}^*$ is positive recurrent under $\theta^{[P]}$. However, this requires some justification. Recall that the server visits the demand points according to the sequence $(1,2,...,d, 1, 2, ..., d, 1, 2, ...)$ under policy $\theta^{[P]}$.  Therefore, at any given time, the server's current node $v$ must lie somewhere on the path between demand points $w$ and $w+1$ (if $w\leq d-1$) or $d$ and $1$ (if $w=d$). Using property (\ref{stab_prop4}), we can assume that the variables $v$ and $w$ associated with state $\vec{u}^*$ do indeed satisfy these constraints, as it is always possible to change the value of $w$ without making any difference to the optimal finite-stage expected cost. Also, if $\vec{u}^*$ is a state with $v\notin D$ (i.e. the server is at an intermediate stage rather than a demand point), it is reasonable to assume that node $v$ is visited during the server's cyclic route under policy $\theta^{[P]}$, since if this is not the case, we can always modify the policy $\theta^{[P]}$ slightly so that node $v$ is visited at some point during the server's route.

Using similar arguments to those given for (AC2), we then have that $\vec{u^*}$ is positive recurrent under $\theta^{[P]}$ and $J_{m,\theta^{[P]}}(\vec{x},\vec{u}^*)\leq J_{\theta^{[P]}}(\vec{x},\vec{u}^*)<\infty$ for any $\vec{x}\in\tilde{S}_m$.  In particular let us consider the state $\vec{z}=(1,1,(0,0,...,0))$. Repeating previous arguments, we have $v_m^{(k)}(\vec{z})\leq J_{\theta^{[P]}}(\vec{z},\vec{u}^*)+ v_m^{(k)}(\vec{u}^*)$. Taking limits as $k\rightarrow\infty$, we obtain $h_m(\vec{u})\geq h_m(\vec{u}^*)\geq -J_{\theta^{[P]}}(\vec{z},\vec{u}^*)$ for all $\vec{u}\in U$. This establishes the required lower bound for states in $U$, and it follows from (\ref{AC3_ineq2}) that the same lower bound also works for all $\vec{x}\in \tilde{S}_m$. Since this argument can be repeated for each $m\in\mathbb{N}$ (establishing a uniform lower bound independent of $m$), we have $-J_{\theta^{[P]}}(\vec{z},\vec{u}^*)\leq \lim\inf_{m\rightarrow\infty} h_m(\vec{x})$ for all $\vec{x}\in\tilde{S}$ as required.

Assumption (AC4) states that $\lim\sup_{m\rightarrow\infty}g_m=:g^*<\infty$ and $g^*\leq g(\vec{x})$ for $\vec{x}\in\tilde{S}$, where $g_m$ is the constant that appears in (\ref{opt_eqs_n}) and the notation $g(\vec{x})$ allows for a possible dependence of the long-run average cost on the initial state $\vec{x}$. By Proposition 8.2.1, Step 3(i) in \cite{sennott1999stochastic}, it is sufficient to show that
\begin{equation}v_m^{(k)}(\vec{x})\leq \lim_{p\rightarrow\infty}v_p^{(k)}(\vec{x})\;\;\;\; \forall\; m,k\in\mathbb{N}_0,\;\vec{x}\in\tilde{S}_m,\label{AC4_ineq}\end{equation}
which follows immediately from property (\ref{stab_prop3}).

Having verified that assumptions (AC1)-(AC4) hold for the modified MDP $\Psi$, we can use the results in \cite{sennott1999stochastic} to conclude that any limit point of a sequence of stationary optimal policies for the finite-state MDPs $(\Psi_0,\Psi_1,\Psi_2,...)$ is optimal for $\Psi$ (Theorem 8.1.1) and furthermore a limit point is guaranteed to exist (Proposition B.5). By the previous arguments, we can compute an optimal policy for any finite-state MDP $\Psi_m$ using value iteration. During the process of value iteration, the functions $v_m^{(k)}(\vec{x})$ for $\vec{x}\in\tilde{S}_m$ are computed using the rule
$$v_m^{(k+1)}(\vec{x})=f(\vec{x})+\min_{a\in A_\vec{x}}\left\{\sum_{\vec{y}\in\tilde{S}_m}p_{\vec{x},\vec{y}}(a)v_m^{(k)}(\vec{y})\right\},\;\;\;\;k\in\mathbb{N}_0.$$
Property (\ref{stab_prop4}) implies that if $\vec{x}$ and $\vec{x}'$ are two states in $\tilde{S}_m$ that differ from each other only in the variable $w$, then any action $a\in A_{\vec{x}}$ that attains the minimum in the equation above for state $\vec{x}$ is also a feasible action that attains the minimum in the corresponding equation for state $\vec{x}'$. Essentially, this means that it is possible to find an optimal stationary policy for $\Psi_m$ that chooses actions independently of the variable $w$. By the previous arguments, the same property also applies to an optimal stationary policy for the infinite-state MDP $\Psi$ (obtained as a limit of the optimal finite-state policies). However, if we have an optimal stationary policy that chooses actions independently of $w$, then the same policy must also be admissible for the MDP formulated in Section \ref{sec:Problem Formulation} with state space $S$. From (AC4) we also know that the long-run average cost under such a policy is finite, implying stability. This completes the proof. \hfill $\Box$\\

\section{DVO Heuristic}\label{AppDVO}

In this appendix we describe the steps of the DVO heuristic as presented in \cite{duenyas1996heuristic}. Note that, where appropriate, we adapt the authors' notation so that it is consistent with the notation used in our paper.
\begin{enumerate}
\item If the server has just finished processing a job at some demand point $i\in D$ then there are two possible cases: either (a) there are still some jobs remaining at $i$ ($x_i>0$) or (b) there are no jobs remaining at $i$ ($x_i=0$). 
\begin{enumerate}[(a)]
\item In the first case, the server can either continue processing jobs at $i$ or switch to some other demand point $j\neq i$. We carry out the following steps:
    \begin{enumerate}
        \item Initialize an empty set $\sigma = \emptyset$.
        \item For each demand point $j$ with $c_j\mu_j\geq c_i\mu_i$, calculate the reward rate $\psi_j$ that would be earned by switching to node $j$, serving it exhaustively, then switching back to node $i$, using
        $$\psi_j = c_j\mu_j\frac{x_j+\lambda_j\delta(i,j)/\tau}{x_j+\mu_j\delta(i,j)/\tau+(\mu_j-\lambda_j)\delta(j,i)/\tau}.$$
        If $\psi_j \geq c_j\mu_j\rho + c_i\mu_i(1-\rho)$, then add $j$ to the set $\sigma$. 
        \item If $\sigma$ is non-empty, then switch to the demand point $j$ with the highest index $\psi_j$ (with ties broken arbitrarily). Otherwise, process one more job at node $i$. 
    \end{enumerate}

\item In the second case, the server can either remain idle at $i$ or switch to some other demand point $j\neq i$. We carry out the following steps:
    \begin{enumerate}
        \item Initialize three empty sets: $\sigma_1 = \emptyset$, $\sigma_2 = \emptyset$ and $\sigma = \emptyset$.
        \item For each demand point $j\neq i$, calculate the reward rate $\phi_j$ in a similar way to part (a) but without including the time taken to switch back from $j$ to $i$, using
        $$\phi_j = c_j\mu_j\frac{x_j+\lambda_j\delta(i,j)/\tau}{x_j+\mu_j\delta(i,j)/\tau}.$$
        If $\phi_j > c_j\mu_j\rho$, then add $j$ to $\sigma_1$; otherwise, add $j$ to $\sigma_2$.
        \item If $\sigma_1$ is non-empty, let $\sigma = \sigma_1$. Otherwise, let $\sigma = \sigma_2$.
        \item Let $j^*$ denote the demand point in $\sigma$ with the highest reward rate $\phi_{j^*}$, with ties broken arbitrarily. If $x_{j^*}>\lambda_{j^*}\delta(j^*,i)/\tau$, then switch to $j^*$. Otherwise, remain idle at $i$. 
    \end{enumerate}

\end{enumerate}
\item If the server has just arrived at a demand point then it immediately begins processing jobs there if there is at least one job waiting. This ensures that, after switching to a new demand point, it must process at least one job there before switching somewhere else. If there are no jobs waiting, then the rule for idling described in step 1(b) is used.
\item If the server is idle at a demand point and a new job arrives in the system, then the rule for idling described in step 1(b) is used.\\
\end{enumerate}

As mentioned in Section \ref{sec:Index heuristic}, there is no need to specify the rule used by the DVO heuristic when the server is at an intermediate stage $i\in N$, since the server is required to continue moving towards a particular demand point (chosen at the previous decision epoch) in this case. Similarly, if the server is at a non-empty demand point $i\in D$ and has \emph{not} just finished processing a job, then this indicates that a service is in progress and the server is required to remain at node $i$.\\ 

\section{Proof of Theorem \ref{pathwise_thm}.}\label{AppB}

\textit{Proof. }
Let $\vec{x}=(v,(x_1,...,x_d))$ be the current state of the system, where $v\in N$. We will use $T_{\textup{arr}}$ to denote the random amount of time until the next job arrives in the system, and $T_{\textup{switch}}$ to denote the amount of time until the server reaches a demand point. According to the rules of the $K$-stop heuristic, a sequence $s$ can only be added to the set $\sigma$ if it satisfies the condition $\frac{\partial}{\partial t}\psi(\vec{x},s,t)\big|_{t=0}\leq 0$. We begin by showing that there always exists at least one sequence that satisfies this condition, and therefore $\sigma_2$ (and, hence, $\sigma$) must be non-empty. Indeed, the set $\mathcal S$ is defined to include all sequences of length $m$, for each $1\leq m\leq K$. Therefore it includes the sequences of length one, $(j)$, for each $j\in D$. Let $s=(j)$, where $j\in D$ is arbitrary. Then, using (\ref{psi_eqn}), we have
\begin{align}\psi(\vec{x},(j),t)&=c_j\mu_j\left\{\frac{T_1(\vec{x},(j),t)}{t+\delta(v,j)/\tau+T_1(\vec{x},(j),t)}\right\}\nonumber\\[10pt]
&=c_j\mu_j\left\{\frac{[x_j+\lambda_j(t+\delta(v,j)/\tau)]/(\mu_j-\lambda_j)}{t+\delta(v,j)/\tau+[x_j+\lambda_j(t+\delta(v,j)/\tau)]/(\mu_j-\lambda_j)}\right\}\nonumber\\[10pt]
&=c_j\mu_j\left\{\frac{x_j+\lambda_j(t+\delta(v,j)/\tau)}{x_j+\mu_j(t+\delta(v,j)/\tau)}\right\}.\nonumber\end{align}
After differentiating, we obtain
\begin{equation}\frac{\partial}{\partial t}\psi(\vec{x},(j),t)=-c_j\mu_j\left\{\frac{x_j(\mu_j-\lambda_j)}{[x_j+\mu_j(t+\delta(v,j)/\tau)]^2}\right\},\label{psi_deriv}\end{equation}
which equals zero if $x_j=0$, and is negative otherwise. Hence, the sequence $s=(j)$ satisfies the condition $\frac{\partial}{\partial t}\psi(\vec{x},s,t)\big|_{t=0}\leq 0$. We can therefore be sure that $\sigma_2$ is non-empty, but in order to proceed we must consider two possible subcases: either (a) the set $\sigma_1$ is non-empty and we choose the sequence in $\sigma_1$ with the highest value of $\psi(\vec{x},s,0)$, or (b) the set $\sigma_1$ is empty, but $\sigma_2$ is non-empty and we choose the sequence in $\sigma_2$ with the highest value of $\psi(\vec{x},s,0)$. In the remainder of this proof we consider these two subcases separately.\\

\textbf{Subcase (a): $\sigma_1$ is non-empty}

For convenience, we will use $\xi(\vec{x},s,t)$ to denote the proportion of time spent processing jobs (as opposed to switching between nodes or idling) while following sequence $s$. That is:
\begin{equation}\xi(\vec{x},s,t):=\dfrac{\sum_{j=1}^{|s|} T_j(\vec{x},s,t)}{t+\sum_{j=1}^{|s|}\left[\delta(s_{j-1},s_j)/\tau+T_j(\vec{x},s,t)\right]},\;\;\;t\geq 0.\label{xi_eqn}\end{equation}
Given that $T_j(\vec{x},s,t)\equiv R_j(\vec{x},s,t)/(c_{s_j}\mu_{s_j})$ for each $j=1,...,|s|$, we can use identical arguments to those in the proof of Lemma \ref{sec3_lem} to show that $\xi(\vec{x},s,t)$ is a monotonic function of $t$. Observe that the condition $\psi(\vec{x},s,0)\geq \gamma(\vec{x},s,0)$ is equivalent to
\begin{equation}\xi(\vec{x},s,0)\geq\rho.\label{xi_condition}\end{equation}
Let $s^*$ denote the sequence in $\sigma_1$ that maximizes $\psi(\vec{x},s,0)$. In this case, according to the rules of the heuristic, the server attempts to take one step along a shortest path to demand point $s_1^*$, where $s_1^*$ is the first element of $s^*$. Let $(i_0,i_1,...,i_k)$ denote a shortest path from $v$ to $s_1^*$, where $k:=\delta(v,s_1^*)$, $i_0:=v$ and $i_k:=s_1^*$. Also let $\vec{x}_j$ denote a state identical to $\vec{x}$ except that the server is located at node $i_j$, for $j=1,2,...,k$.   It is useful to note that
\begin{equation}\xi(\vec{x}_j,s^*,0)=\xi(\vec{x}_k,s^*,(k-j)/\tau)\;\;\;\;\forall j\in\{1,2,...,k-1\}.\label{proof_eqn1}\end{equation}
This is because $\xi(\vec{x}_j,s^*,0)$ represents the proportion of time spent processing jobs given that the server begins at node $i_j$ and immediately begins traveling along a path of length $(k-j)$ in order to reach node $s_1^*$, while $\xi(\vec{x}_k,s^*,(k-j)/\tau)$ is the corresponding proportion given that the server begins at node $i_k=s_1^*$ but waits for $(k-j)/\tau$ time units before beginning to process jobs there. In terms of the proportion of time spent processing jobs while following sequence $s^*$, these two quantities are the same. In the case $|s^*|=1$ (that is, $s^*=(j)$ for some $j\in D$) we infer from (\ref{psi_deriv}) that
\begin{equation}\frac{\partial}{\partial t}\xi(\vec{x},(j),t)=\frac{1}{c_j\mu_j}\frac{\partial}{\partial t}\psi(\vec{x},(j),t)\leq 0,\label{xi_deriv}\end{equation}
and hence, using (\ref{proof_eqn1}), we can be assured that the condition $\xi(\vec{x},s^*,0)\geq \rho$ (equivalent to $\psi(\vec{x},s^*,0)\geq \gamma(\vec{x},s^*,0)$) remains satisfied as the server moves towards $s_1^*$. Therefore, according to the rules of the heuristic, $s^*$ remains included in $\sigma_1$ at all stages while the server moves from $v$ to $s_1^*$. On the other hand, consider the case $|s^*|\geq 2$. In this case the rules of the heuristic imply that the extra condition $\xi(\vec{y},s^*,0)\geq\rho$ must be satisfied, where $\vec{y}$ is equivalent to $\vec{x}_k$ in the notation of this proof. By the same reasoning used to derive (\ref{proof_eqn1}), we have
\begin{equation}\xi(\vec{x},s^*,0)=\xi(\vec{x}_k,s^*,k/\tau).\label{proof_eqn2}\end{equation}
Since $\xi(\vec{x}_k,s^*,t)$ is a monotonic function of $t$, we can infer from (\ref{proof_eqn1}) and (\ref{proof_eqn2}) that
\begin{equation}\min\left\{\xi(\vec{x},s^*,0),\;\xi(\vec{x}_k,s^*,0)\right\}\leq \xi(\vec{x}_j,s^*,0)\leq \max\left\{\xi(\vec{x},s^*,0),\;\xi(\vec{x}_k,s^*,0)\right\}\;\;\;\;\forall j\in\{1,2,...,k\}.\label{proof_eqn3}\end{equation}
Given that $s^*\in\sigma_1$, the rules of the heuristic imply that $\xi(\vec{x},s^*,0)\geq \rho$ and $\xi(\vec{x}_k,s^*,0)\geq \rho$, so from (\ref{proof_eqn3}) it follows that $\xi(\vec{x}_j,s^*,0)\geq\rho$ for each $j=1,2,...,k$. Also, given that $s^*\in\sigma_1$, the heuristic rules imply $\frac{\partial}{\partial t}\psi(\vec{x},s^*,t)\big|_{t=0}\leq 0$. Applying the same reasoning to the function $\psi$ that we used for $\xi$, we have
\begin{align}&\psi(\vec{x}_j,s^*,0)=\psi(\vec{x}_k,s^*,(k-j)/\tau)\;\;\;\;\forall j\in\{1,2,...,k\},\label{proof_eqn_3a}\\
&\psi(\vec{x},s^*,0)=\psi(\vec{x}_k,s^*,k/\tau).\label{proof_eqn_3b}\end{align}
From the proof of Lemma \ref{sec3_lem} we know that $\frac{\partial}{\partial t}\psi(\vec{x},s^*,t)$ has the same sign for all $t\geq 0$. Hence, from (\ref{proof_eqn_3b}) it follows that $\frac{\partial}{\partial t}\psi(\vec{x}_k,s^*,t)\big|_{t=0}\leq 0$ and then from (\ref{proof_eqn_3a}) it follows that $\frac{\partial}{\partial t}\psi(\vec{x}_j,s^*,t)\big|_{t=0}\leq 0$ for each $j=1,2,...,k$. By collating the arguments given thus far, we can say that as the server travels from $v$ to $s_1^*$ along the shortest path $(i_0,i_1,...,i_k)$, the sequence $s^*$ always satisfies the necessary conditions to be included in $\sigma_1$. 

Now suppose that the server arrives at node $i_1$ before time $T_{\textup{arr}}$. If $i_1=s_1^*$ (i.e. the shortest path from $v$ to $s_1^*$ is of length 1), then there is nothing further to prove, as the server has completed its journey to demand point $s_1^*$ via the shortest possible path. On the other hand, if $i_1\neq s_1^*$ (i.e. the shortest path from $v$ to $s_1^*$ has length greater than 1), then we must continue. Given that $\frac{\partial}{\partial t}\psi(\vec{x},s^*,t)\big|_{t=0}\leq 0$ and using Lemma \ref{sec3_lem} again, it must be the case that $\psi(\vec{x},s^*,0)\leq \psi(\vec{x}_1,s^*,0)$; in other words, the average reward does not decrease after the server takes a step along the shortest path to $s_1^*$. By the previous arguments, $s^*$ is still included in $\sigma_1$ when the server is located at $i_1$. However, we cannot be sure that $s^*$ is chosen by the heuristic at state $\vec{x}_1$. 

If $s^*$ is chosen by the heuristic at $\vec{x}_1$, then the server continues attempting to move along the shortest path to node $s_1^*$. On the other hand, suppose the heuristic chooses an alternative sequence $\tilde{s}\in\sigma_1$ at $\vec{x}_1$, and let $\tilde{s}_1$ be the first demand point in sequence $\tilde{s}$. If $\tilde{s}_1=s_1^*$, then (once again) we have no difficulties, as the server continues attempting to move along the shortest path to $s_1^*$. In the rest of this subcase we assume a non-trivial case where $\tilde{s}_1\neq s_1^*$.

We can show that although $\tilde{s}_1\neq s_1^*$, the server's movement from node $v$ to $i_1$ still qualifies as a step along a shortest path from $v$ to $\tilde{s}_1$. In other words, even if the server now prefers to move towards a different demand point, the step from $v$ to $i_1$ was still a step in the right direction. To see this, first note that 
\begin{equation}\psi(\vec{x},\tilde{s},0)\leq \psi(\vec{x},s^*,0)\leq \psi(\vec{x}_1,s^*,0)\leq \psi(\vec{x}_1,\tilde{s},0),\label{proof_eqn3b}\end{equation}
where the first inequality is due to the fact that $s^*$ is preferred to $\tilde{s}$ at state $\vec{x}$, the second inequality is due to the fact that $\frac{\partial}{\partial t}\psi(\vec{x},s^*,t)\big|_{t=0}\leq 0$ (as stated earlier) and the third inequality is due to the fact that $\tilde{s}$ is preferred to $s^*$ under the new state $\vec{x}_1$. The rules of the heuristic state that if the average rewards for two sequences are equal, then a sequence is chosen according to a fixed priority ordering of the demand points in $D$. Therefore either the first inequality or the third inequality in (\ref{proof_eqn3b}) must be strict, as if they both hold with equality then the same sequence (either $s^*$ or $\tilde{s}$) must be chosen at both $\vec{x}$ and $\vec{x}_1$. We conclude that 
\begin{equation}\psi(\vec{x},\tilde{s},0)<\psi(\vec{x}_1,\tilde{s},0),\label{proof_eqn3c}\end{equation}
implying that the average reward for sequence $\tilde{s}$ has increased after moving from $v$ to $i_1$. Let $\tilde{\vec{x}}$ denote a state identical to $\vec{x}$ except that the server is located at node $\tilde{s}_1$ instead of $v$. By definition, we have
\begin{align}&\psi(\vec{x},\tilde{s},0)=\psi(\tilde{\vec{x}},\tilde{s},\delta(v,\tilde{s}_1)/\tau),\label{proof_eqn5}\\
&\psi(\vec{x}_1,\tilde{s},0)=\psi(\tilde{\vec{x}},\tilde{s},\delta(i_1,\tilde{s}_1)/\tau).\label{proof_eqn6}\end{align}
Suppose (for a contradiction) that $\delta(i_1,\tilde{s}_1)\geq\delta(v,\tilde{s}_1)$. If these two distances are equal then from (\ref{proof_eqn5})-(\ref{proof_eqn6}) we have $\psi(\vec{x},\tilde{s},0)=\psi(\vec{x}_1,\tilde{s},0)$, giving a contradiction with (\ref{proof_eqn3c}). On the other hand, if the inequality is strict (that is, $\delta(i_1,\tilde{s}_1)>\delta(v,\tilde{s}_1)$) then using Lemma \ref{sec3_lem} we infer that $\frac{\partial}{\partial t}\psi(\vec{x}_1,\tilde{s},t)\big|_{t=0}>0$. Hence, according to the rules of the heuristic, the sequence $\tilde{s}$ could not have been chosen under state $\vec{x}_1$. We conclude that $\delta(i_1,\tilde{s}_1)< \delta(v,\tilde{s}_1)$ and therefore the server's step from $v$ to $i_1$ represents a step along a shortest path from $v$ to $\tilde{s}_1$.\\

\textbf{Subcase (b): $\sigma_1$ is empty}

From the arguments at the beginning of the proof we know that even if $\sigma_1$ is empty, $\sigma_2$ must be non-empty. As in case (a), let $s^*$ denote the sequence chosen by the heuristic under state $\vec{x}$, and let $(i_0,i_1,...,i_k)$ denote a shortest path from $v$ to $s_1^*$, where $k=\delta(v,s_1^*)$, $i_0=v$ and $i_k=s_1^*$. Also let $\vec{x}_j$ denote a state identical to $\vec{x}$ except that the server is located at node $i_j$, for $j=1,2,...,k$. After the server moves from node $v$ to $i_1$, there are two possible scenarios: either $\sigma_1$ remains empty, or it becomes non-empty. In the first scenario, we can simply repeat the relevant arguments used in subcase (a) to show that, even if the heuristic chooses some alternative sequence $\tilde{s}\in\sigma_2$ under state $\vec{x}_1$, it must be the case that $\delta(i_1,\tilde{s}_1)<\delta(v,\tilde{s}_1)$, and therefore the server's movement from $v$ to $i_1$ qualifies as a step along a shortest path from $v$ to $\tilde{s}_1$. This is due to the fact that $\frac{\partial}{\partial t}\psi(\vec{x}_1,\tilde{s},t)\big|_{t=0}$ must be non-positive in order for $\tilde{s}$ to be chosen. In the rest of this part we consider the second scenario, where $\sigma_1$ becomes non-empty.

Suppose $\tilde{s}\in\sigma_1$ and the heuristic chooses sequence $\tilde{s}$ under state $\vec{x}_1$. It may be the case that $\tilde{s}_1=s_1^*$, in which case the server simply keeps following the same path. However, if $\tilde{s}_1\neq s_1^*$, we can again show that $\delta(i_1,\tilde{s}_1)<\delta(v,\tilde{s}_1)$. To see this, let $\tilde{\vec{x}}$ denote the state identical to $\vec{x}$ except that the server is located at node $\tilde{s}_1$. Due to the rules of the heuristic, we must have
\begin{equation}\xi(\vec{x}_1,\tilde{s},0)\geq \rho\label{subcase_b_eqn}\end{equation}
and additionally, if $|\tilde{s}|\geq 2$, then
\begin{equation}\xi(\tilde{\vec{x}},\tilde{s},0)\geq \rho.\label{subcase_b_eqn2}\end{equation}
Note that $\xi(\vec{x}_1,\tilde{s},0)=\xi(\tilde{\vec{x}},\tilde{s},\delta(i_1,\tilde{s}_1)/\tau)$. Suppose (for a contradiction) that $\delta(i_1,\tilde{s}_1)\geq \delta(v,\tilde{s}_1)$. Given that sequence $\tilde{s}$ was not included in $\sigma_1$ when the server was at state $\vec{x}$, at least one of the conditions $\frac{\partial}{\partial t}\psi(\vec{x},\tilde{s},t)\big|_{t=0}\leq 0$, $\xi(\vec{x},\tilde{s},0)\geq\rho$ and $\xi(\tilde{\vec{x}},\tilde{s},0)\geq\rho$ (where the latter only applies if $|\tilde{s}|\geq 2$) must fail to hold. However, given that $\tilde{s}$ is chosen at $\vec{x}_1$ and $\psi(\vec{x}_1,\tilde{s},0)\equiv\psi(\vec{x},\tilde{s},1/\tau)$, it must be the case (using Lemma \ref{sec3_lem}) that $\frac{\partial}{\partial t}\psi(\vec{x}_1,\tilde{s},t)\big|_{t=0}$ has the same sign as $\frac{\partial}{\partial t}\psi(\vec{x},\tilde{s},t)\big|_{t=0}$, so the derivative condition holds. We also have $\xi(\tilde{\vec{x}},\tilde{s},0)\geq\rho$ when $|\tilde{s}|\geq 2$ from (\ref{subcase_b_eqn2}), so we can proceed to assume that $\xi(\vec{x},\tilde{s},0)<\rho$, which is equivalent to $\xi(\tilde{\vec{x}},\tilde{s},\delta(v,\tilde{s}_1)/\tau)<\rho$. Hence, we have
\begin{equation}\xi(\tilde{\vec{x}},\tilde{s},\delta(v,\tilde{s}_1)/\tau)<\rho\leq \xi(\tilde{\vec{x}},\tilde{s},\delta(i_1,\tilde{s}_1)/\tau)\label{subcase_b_eqn4}\end{equation}
but also (by assumption)
\begin{equation}0<\delta(v,\tilde{s}_1)\leq \delta(i_1,\tilde{s}_1),\label{subcase_b_eqn3}\end{equation}
implying that $\xi(\tilde{\vec{x}},\tilde{s},t)$ is increasing with $t$. If $|\tilde{s}|=1$ then this gives a contradiction, since it was shown in (\ref{xi_deriv}) that $\xi(\cdot)$ is non-increasing with $t$ for sequences of length one. On the other hand, if $|\tilde{s}|\geq 2$, we modify the right-hand side of (\ref{subcase_b_eqn4}) and combine (\ref{subcase_b_eqn})-(\ref{subcase_b_eqn2}) to give
$$\xi(\tilde{\vec{x}},\tilde{s},\delta(v,\tilde{s}_1)/\tau)<\rho\leq \min\{\xi(\tilde{\vec{x}},\tilde{s},0),\;\xi(\tilde{\vec{x}},\tilde{s},\delta(i_1,\tilde{s}_1)/\tau)\},$$
which, along with (\ref{subcase_b_eqn3}), contradicts the fact that $\xi(\tilde{\vec{x}},\tilde{s},t)$ is a monotonic function of $t$. We conclude that $\delta(i_1,\tilde{s}_1)<\delta(v,\tilde{s}_1)$ as required.\\

We can repeat the arguments given in subcases (a) and (b) to show that any movement of the server from one intermediate stage to another (prior to $\min\{T_{\textup{arr}},\;T_{\textup{switch}}\}$) qualifies as a step along a shortest path from $v$ to some particular demand point $j^*$. The key point is that even though the sequences chosen at the various intermediate stages may not always begin with node $j^*$, the distance from the currently-occupied node to $j^*$ is always reduced after each step, so the complete path is indeed a shortest path from $v$ to $j^*$. Given that the number of nodes in the network is finite, the demand point $j^*$ must eventually be reached if no new jobs arrive in the meantime. This completes the proof. \hfill $\Box$\\

\section{Proof of Corollary \ref{pathwise_corollary}.}\label{AppC}

\textit{Proof. }We assume that the server is initially located at some intermediate stage $v\in N$ and use $M=\max_{\{i\in N,\;j\in D\}}\delta(i,j)$ to denote the maximum distance between an intermediate stage and a demand point. Let $\alpha$ denote the number of new jobs that arrive in the system before the server reaches a demand point, given that the $K$-stop heuristic is followed. By conditioning on $\alpha$, we have
\begin{equation}\mathbb{E}[T_{\textup{switch}}]=\sum_{k=0}^\infty \mathbb{E}[T_{\textup{switch}}\;|\;\alpha=k]\;\mathbb{P}(\alpha=k).\label{cor_eq1}\end{equation}
Due to Theorem \ref{pathwise_thm} we know that if the server is at an intermediate stage, then until the next new job arrives it attempts to move along a shortest path to some particular demand point $j^*$. Since (prior to $T_{\textup{switch}}$) the server is always attempting to move, the expected amount of time until the next event (either a switch or the arrival of a new job in the system) is always $(\Lambda+\tau)^{-1}$, where $\Lambda=\sum_{i\in D}\lambda_i$. In the worst case, the number of nodes that must be traversed in order to reach $j^*$ is $M$. Hence, we can form an upper bound for $\mathbb{E}[T_{\textup{switch}}\;|\;\alpha=0]$:
$$\mathbb{E}[T_{\textup{switch}}\;|\;\alpha=0]\leq\frac{M}{\Lambda+\tau}.$$
Extending this argument, we can obtain an upper bound for $\mathbb{E}[T_{\textup{switch}}\;|\;\alpha=k]$ (for $k\geq 1$) by supposing that every time a new job arrives in the system, the server changes direction and attempts to move to a demand point $M$ nodes away, and manages to complete $(M-1)$ of these switches before a new job arrives in the system and forces it to change direction again. Suppose this pattern continues until $k$ arrivals have occurred, at which point it manages to complete $M$ switches without interruption and reaches a demand point. In this scenario, the total number of switches made is $k(M-1)+M$ and the total number of new jobs arriving is $k$, so the total number of system events is $(k+1)M$. Hence:
$$\mathbb{E}[T_{\textup{switch}}\;|\;\alpha=k]\leq \frac{(k+1)M}{\Lambda+\tau},\;\;\;\;k\geq 0.$$
Next, consider the probabilities $\mathbb{P}(\alpha=k)$. Each time a system event occurs (prior to $T_{\textup{switch}}$), there is a probability of $\tau/(\Lambda+\tau)$ that this is a switch rather than the arrival of a new job. We can obtain an upper bound for $\mathbb{P}(\alpha=0)$ by supposing that the server only needs to complete one switch in order to reach its intended demand point. Hence:
$$\mathbb{P}(\alpha=0)\leq\frac{\tau}{\Lambda+\tau}.$$
On the other hand, the largest possible probability of an arrival occurring before the server reaches a demand point is $1-(\tau/(\Lambda+\tau))^M$, since this represents the case where the server must complete $M$ switches in order to reach its intended demand point. Putting these arguments together, we have the following upper bound:
$$\mathbb{P}(\alpha=k)\leq\left[1-\left(\frac{\tau}{\Lambda+\tau}\right)^M\right]^k \left(\frac{\tau}{\Lambda+\tau}\right),\;\;\;\;k \geq 0.$$
Let $p:=1-(\tau/(\Lambda+\tau))^M$ for notational convenience. Then, using (\ref{cor_eq1}), we have
\begin{align}\mathbb{E}[T_{\textup{switch}}]&\leq\sum_{k=0}^\infty \frac{(k+1)M}{\Lambda+\tau}\; p^k \left(\frac{\tau}{\Lambda+\tau}\right)\nonumber\\[8pt]
&=\frac{\tau M}{(\Lambda+\tau)^2}\sum_{k=0}^\infty (k+1)p^k\nonumber\\[8pt]
&=\frac{\tau M}{(\Lambda+\tau)^2}\cdot\frac{1}{(1-p)^2}\nonumber\\[8pt]
&=\frac{\tau M}{(\Lambda+\tau)^2}\left(\frac{\Lambda+\tau}{\tau}\right)^{2M}\nonumber\\[8pt]
&=\frac{M}{\tau}\left(\frac{\Lambda+\tau}{\tau}\right)^{2(M-1)}.\end{align}
This completes the proof. We also note that the bound holds with equality if and only if $M=1$.  \hfill $\Box$\\

\black{\section{Proof of Theorem \ref{homogeneous_thm}.}\label{App_homsystem}}

\black{\textit{Proof. }We will use $\lambda$, $\mu$ and $c$ to denote the common arrival rate, service rate and holding cost (respectively) for all job types. Firstly, let $\vec{x}$ be a state under which the server is located at a demand point $i\in D$ with $x_i\geq 1$. Recall from step 1(a) of the $K$-stop heuristic algorithm that $\mathcal S$ is the set of sequences of the form $s=(s_1,s_2,...,s_m)$, where $1\leq m\leq K$, $s_j\in D$ for each $j\in\{1,2,...,m\}$, $s_1\neq v$ and $s_i\neq s_j$ for any pair of elements $s_i,s_j\in s$ with $i\neq j$. For any sequence $s\in\mathcal S$, $t\geq 0$ and $j=\in\{1,...,|s|\}$ it is clear from (\ref{phi_eqn}) that
$$\phi_j(\vec{x},s,t)<\dfrac{\sum_{k=1}^j R_k(\vec{x},s,t)}{\sum_{k=1}^j T_k(\vec{x},s,t)}=\dfrac{c\mu\sum_{k=1}^j T_k(\vec{x},s,t)}{\sum_{k=1}^j T_k(\vec{x},s,t)}=c\mu.$$
On the other hand, the quantity $\beta_j(\vec{x},s,t)$ in (\ref{beta_eqn}) is equal to $c\mu$. Therefore the condition $\phi_j(\vec{x},s,t)\geq\beta(\vec{x},s,t)$ is not satisfied for any $s\in\mathcal S$, and according to the rules of the heuristic the server should remain at node $i$. The intuitive explanation for this is that the server earns rewards at the maximum possible rate $c\mu$ by remaining at its current node, so there is no reason to switch to another node. It follows that, under the $K$-stop policy, the server always remains at a demand point until all jobs have been processed.}

\black{We can also show that the server visits all demand points infinitely often under the $K$-stop policy. Indeed, suppose (for a contradiction) that there exists some non-empty subset $D_0\subset D$ such that demand points in $D_0$ are visited only finitely many times. We also define $D_1:=D\setminus D_0$ as the subset of demand points that are visited infinitely often. The subset $D_1$ must be non-empty because, as noted in the proof of Theorem \ref{pathwise_thm}, sequences $s$ of length one always satisfy the condition $\frac{\partial}{\partial t}\psi(\vec{x},s,t)\big|_{t=0}\leq 0$, and therefore the set $\sigma$ constructed in step 2 of the heuristic algorithm must be non-empty, which implies that if the server is at an intermediate stage then it will move towards a demand point. 
Each time the server selects a demand point in $i\in D_1$ to move to, there is a positive probability that no new jobs arrive at any of the other demand points $j\in D_1\setminus\{i\}$ while it switches to $i$ and processes jobs there. Therefore the system must eventually reach a state in which there are no jobs at any of the demand points in $D_1$. Let $\vec{x}$ be a state with $x_j=0$ for all $j\in D_1$ and let $\mathcal S_1$ denote the set of all sequences in $\mathcal S$ that involve visiting only demand points in $D_1$. Also, define }
\black{$$\psi_1^{\max}:=\max_{s\in\mathcal S_1} \psi(\vec{x},s,0).$$}
\black{It is clear from the definition of $\psi(\vec{x},s,0)$ in (\ref{psi_eqn}) that $\psi_1^{\max}$ is strictly smaller than $c\mu$. Next, consider an arbitrary demand point $j\in D_0$ and consider the sequence of length one, $s=(j)$. We have assumed that demand point $j$ is visited only finitely many times and therefore, in the long run, $x_j$ tends to infinity. From (\ref{psi_eqn}) it follows that the index $\psi(\vec{x},(j),0)$ tends towards $R_1(\vec{x},(j),0)/T_1(\vec{x},(j),0)=c\mu$. Hence, at some point $\psi(\vec{x},(j),0)$ must exceed $\psi_1^{\max}$. Similar reasoning can be used to show that any sequence $s$ that begins by visiting a demand point in $D_1$ must eventually become inferior (in terms of index value) to the sequence $(j)$, where $j\in D_0$. The rules of the heuristic then imply that the server selects a sequence which begins by visiting one of the demand points in $D_0$, which yields a contradiction. We conclude that all demand points must be visited infinitely often under the $K$-stop policy. Standard arguments based on stability in polling systems (see the comments following the statement of Theorem \ref{stability_theorem}) then imply that the system is stable, which proves the first statement in the theorem.} 

\black{Next, we make the additional assumption that $V$ is a complete graph, which implies that the server can move between any two demand points $i,j\in D$ by traversing a single edge of the network, without having to pass through any intermediate stages. We will show that, with this extra assumption (in addition to the homogeneity assumptions made already), the $K$-stop policy causes the server to act in the following way:}
\black{\begin{enumerate}[(i)]
\item If the server is at a non-empty demand point $i\in D$, then it remains there.
\item If the server is at an empty demand point $i\in D$, then it switches to the demand point $j\neq i$ with the largest number of jobs present. (If there is a tie, an arbitrary choice of $j$ can be made.)
\end{enumerate}}
\black{Note that the rules (i)-(ii) are sufficient to completely specify the actions chosen by the $K$-stop policy, since it cannot ever visit an intermediate stage. Property (i) was already established in the first part of the proof (the additional assumption that $V$ is a complete graph does not alter the argument given previously). To establish property (ii), we again recall the previous argument and note that the set $\sigma$ constructed in step 2 of the heuristic algorithm is non-empty, because it includes sequences of length one. It follows that, given some state $\vec{x}$ under which the server is at an empty demand point, it selects the sequence $s$ that maximizes $\psi(\vec{x},s,0)$. In general, the sequence selected by the $K$-stop policy could be of any length between $1$ and $K$. However, we can show that if a sequence of length greater than one is selected, then the demand points in the sequence must be visited in descending order of $x_i$; that is, the server prioritizes the demand point with the greatest number of jobs. To see this, note that for a sequence $s$ of length $m$ (where $1\leq m\leq K$), the index $\psi(\vec{x},s,0)$ can be expressed as}
\black{\begin{equation}\psi(\vec{x},s,0)=\frac{cu\sum_{k=1}^m T_k(\vec{x},s,0)}{m/\tau+\sum_{k=1}^m T_k(\vec{x},s,0)},\label{psi_proof_eqn}\end{equation}}
\black{where the index $k$ corresponds to the position of demand point $s_k$ in the sequence, so $T_1(\vec{x},s,0)$ is the amount of time spent at the first demand point in the sequence (under the fluid model), etc. Suppose we have a sequence in which the demand points are \emph{not} ordered in descending order of $x_i$. This implies that there must be some $k\in\{1,...,m-1\}$ such that $x_{s_k}<x_{s_{k+1}}$. We will show that the index in (\ref{psi_proof_eqn}) would increase if we swapped the positions of demand points $s_k$ and $s_{k+1}$ in the sequence. Indeed, let $y=x_{s_{k+1}}-x_{s_k}>0$. After performing the swap, $T_k(\vec{x},s,0)$ increases by $y/(\mu-\lambda)$ since, under the fluid model, the amount of time taken to process a single job is $1/(\mu-\lambda)$. On the other hand, $T_{k+1}(\vec{x},s,0)$ increases by $\lambda y/((\mu-\lambda)^2)-y/(\mu-\lambda)$, where the first term $\lambda y/((\mu-\lambda)^2)$ is due to the extra arrivals that occur during the extra $y/(\mu-\lambda)$ time units spent at node $s_k$ and the second term $-y/(\mu-\lambda)$ is due to the fact that there are $y$ fewer jobs at $s_{k+1}$ after performing the swap.} 
\black{Hence, the overall increase in $T_k(\vec{x},s,0)+T_{k+1}(\vec{x},s,0)$ is $\lambda y/((\mu-\lambda)^2)$, which is non-negative. Furthermore, $T_{l}(\vec{x},s,0)$ increases for all $l\geq k+2$ following the swap, due to the fact that $T_k(\vec{x},s,0)+T_{k+1}(\vec{x},s,0)$ increases and therefore extra arrivals occur at these nodes before the server arrives. The result is that the sum $\sum_{k=1}^m T_k(\vec{x},s,0)$ increases following the swap, and therefore the index in (\ref{psi_proof_eqn}) also increases. By extending this argument, we can claim that if $s$ is a sequence in which the demand points are \emph{not} visited in descending order of $x_i$, then it is beneficial to perform a sequence of swaps until they are in descending order. It follows that the sequence that maximizes (\ref{psi_proof_eqn}) must be a sequence in which the first node visited has the greatest number of jobs.}

\black{Having established that properties (i) and (ii) hold under the $K$-stop policy, the next task is to show that these properties also hold under an optimal policy, and thus the $K$-stop policy is optimal. We can show this using a sample path argument. To establish (i), suppose the system is in state $\vec{x}$ with $v(\vec{x})=i\in D$, $x_i\geq 1$. Consider two continuous-time processes, referred to as $P_1$ and $P_2$ for convenience, both initialized in state $\vec{x}$. In the continuous-time setting, an action $a(t)\in D$ must be specified by the server for each time point $t\in\mathbb{R}_{\geq 0}$. In $P_1$ we suppose that the server follows a policy whereby the action at $t=0$ is to switch to demand point $j\neq i$. In $P_2$, on the other hand, the server remains at node $i$ until a service completion occurs there, then copies the same sequence of actions that $P_1$ chose starting from $t=0$, except that it does not choose $i$ again (skipping these actions as necessary) until $P_1$ has completed its first service at $i$. Using standard stochastic coupling arguments, we assume that new job arrivals occur at the same times under both processes, and furthermore the times required for service completions and switch completions are the same under both processes, although (due to the difference in server behaviors) a particular service or switch that begins at time $t\geq 0$ in $P_1$ may begin at a different time $t'\neq t$ under $P_2$. As an example, suppose the trajectory of (action, duration, completion) pairs under $P_1$ evolves as follows:}
\black{$$(j^{\text{switch}},T_0,\checkmark),(j^{\text{serv}},T_1,\times),(i^{\text{switch}},T_2,\checkmark),(i^{\text{serv}},T_3,\times),(k^{\text{switch}},T_4,\checkmark),...,(i^{\text{serv}},T_m,\checkmark),...$$}
\black{which indicates that the server begins by switching to $j$ for $T_0$ time units and this action is completed (denoted by `$\checkmark$'), then it serves a job at node $j$ for $T_1$ time units but this service is not completed (denoted by `$\times$') because it interrupts the service to switch to node $i$, and the switch to $i$ takes $T_2$ units and is successfully completed, etc. We use $m$ to represent the sequence position of the first completed service at $i$. Then the corresponding trajectory under $P_2$ would begin with}
\black{$$(i^{\text{serv}},T_m,\checkmark),(j^{\text{switch}},T_0,\checkmark),(j^{\text{serv}},T_1,\times),(k^{\text{switch}},T_4,\checkmark),...$$}
\black{Let $\hat{T}$ denote the total amount of time that the server spends switching to demand point $i$ (whether the switches are completed or not) in $P_1$ before it eventually completes a service there. In $P_2$ we assume that, after the server has finished copying the non-$i$ actions of the server in $P_1$ prior to $T_0+...+T_{m-1}$, it then chooses action $i$ for a further $\hat{T}$ time units. This ensures that both processes are in the same state at time $T_0+...+T_m$, and they continue to evolve identically from that point onwards.}

\black{We can compare the total costs incurred by $P_1$ and $P_2$ under the above coupling construction. We note that:}
\black{\begin{enumerate}[(a)]
\item Any service completion that occurs in $P_1$ at time $t\geq 0$, except for the first job at demand point $i$, occurs with a delay of no more than $T_m$ time units in $P_2$. Therefore $P_2$ incurs extra holding costs of (at most) $c\times R\times T_m$ due to the delayed processing of these jobs, where $R$ is the total number of jobs processed in $P_1$ before processing the first job at node $i$.
\item The first job at demand point $i$ is processed at time $T_m$ in $P_2$, but is not processed until time $T_0+...+T_m$ in $P_1$. Therefore $P_1$ incurs an extra holding cost of $c\times (T_0+...+T_{m-1})$ due to the delayed processing of this job. 
\end{enumerate}}

\black{Note that the set $\{T_0,...,T_{m-1}\}$ includes the times needed to process $R$ jobs at other demand points before the first job is processed at $i$ under $P_1$. Hence, given that the service times of all jobs are independent and identically distributed, we have $\mathbb{E}[T_0+...+T_{m-1}]\geq \mathbb{E}[R\times T_m]$, so the extra costs incurred in $P_1$ are indeed higher than those in $P_2$. Therefore, switching away from a non-empty demand point is a suboptimal action in the context of minimizing total costs over a long time horizon. Intuitively, the suboptimality occurs because the server is forced to perform extra switching actions in $P_1$.}

\black{Next, we move on to property (ii) and aim to show that if the server is at an empty demand point then it is optimal to switch to the demand point with the largest number of jobs. It is sufficient to show that for any state $\vec{x}=(i,(x_1,...,x_d))$ with $x_i<x_j$ for some $j\in D$, the minimal expected total cost (over a long time horizon) starting from state $\vec{x}$ is greater than it would be if $x_i$ and $x_j$ were swapped. 
We will again use a sample path argument. Let $P_1$ and $P_2$ be two continuous-time processes initialized at in states $\vec{x}=(i,(x_1,...,x_d))$ and $\vec{x}'=(i,(x_1',...,x_d'))$ respectively, where $\vec{x}'$ is identical to $\vec{x}$ except that $x_i'=x_j$ and $x_j'=x_i$, and we assume $x_i<x_j$. 
The two processes are coupled in the same way as in (i), so that new job arrivals occur at the same times under both processes and the times needed for successful service and switch completions are also equivalent.}

\black{Suppose that in $P_1$ the server remains at node $i$ until $y\geq 0$ jobs have been processed before attempting to switch to a different demand point. 
In $P_2$, the server begins by remaining at node $i$ until $y+1$ services have been completed before switching to another demand point, so we force it to process an extra job before switching. It then copies the same actions as the server in $P_1$, but `skips' any actions $j$ chosen in $P_1$ (recall that $P_1$ has at least one extra job at $j$) until eventually $P_1$ completes its first service at $j$. For example, if $y=1$ then the trajectories under $P_1$ and $P_2$ could be as follows:}
\black{\begin{align}&P_1:\;(i^{\text{serv}},T_0,\checkmark),(j^{\text{switch}},T_1,\times),(k^{\text{switch}},T_2,\checkmark),...,(j^{\text{serv}},T_m,\checkmark),...\nonumber\\
&P_2:\;(i^{\text{serv}},T_0,\checkmark),(i^{\text{serv}},T_m,\checkmark),(k^{\text{switch}},T_2,\checkmark),...\nonumber\end{align}}
\black{for some $m\geq 0$. The important thing in this case is to make a comparison between the first job processed at $j$ in $P_1$ and the $(y+1)^{\text{th}}$ job processed at $i$ in $P_2$, and our coupling construction assumes these jobs have the same service time (we can think of them as the same job, which begins at $j$ in $P_1$ but begins at $i$ in $P_2$). Similarly to (i), we also let $\hat{T}$ denote the total amount of time spent by the server switching to demand point $j$ in $P_1$ before it completes a service there and specify that in $P_2$, the server should spend $\hat{T}$ time units choosing action $j$ after it finishes copying the non-$j$ actions of the server in $P_1$ prior to $T_0+...+T_{m-1}$. This ensures that both processes are in the same state at time $T_0+...+T_m$.}

\black{The comparison between the costs incurred by the two processes then works in a similar way to the comparison in (i). In $P_2$, there is a certain number of jobs that have their services delayed by (at most) $T_m$ time units compared to $P_1$, but the extra holding costs for these are collectively smaller than the extra holding cost incurred by $P_1$ as a result of processing the first job at $j$ later than the processing of the $(y+1)^{\text{th}}$ job at $i$ in $P_2$.}

\black{We note that the argument given above becomes most meaningful in the case where the server in $P_1$ remains at the initial node $i$ until it is empty (which it must do under an optimal policy, according to our previous argument). In this case the server in $P_2$ can still process an extra job at $i$, but the server in $P_1$ has to switch at least once before processing the corresponding extra job at node $j$. We conclude that it is always beneficial for the server to be at the demand point with the largest number of jobs present, and it follows that if the server is at an empty demand point then it should attempt to switch to the demand point with the largest number of jobs. We have shown that, in all cases, the actions of the $K$-stop policy are identical to those selected by an optimal policy. \hfill $\Box$}\\

\section{Methods for generating the parameters for the numerical experiments in Section \ref{sec:numerical}}\label{AppD}
For each of the 10,000 instances considered in Section \ref{sec41}, the system parameters are randomly generated as follows:
\begin{itemize}
    \item The size of the left-hand cluster, $d_1$, is sampled unbiasedly from the set $\{1,2,3,4\}$.
    \item The size of the right-hand cluster, $d_2$, is sampled unbiasedly from the set $\{1,2,3,4\}$.
    \item The number of intermediate stages, $n$, is sampled unbiasedly from the set $\{1,2,3,4,5,6\}$. 
    \item The overall traffic intensity, $\rho$, is sampled from a continuous uniform distribution between 0.1 and 0.9. Subsequently, the job arrival rates $\lambda_i$ and processing rates $\mu_i$ for $i\in D$ are generated as follows:
    \begin{itemize}
        \item Each processing rate $\mu_i$ is initially sampled from a continuous uniform distribution between 0.1 and 0.9.
        \item For each demand point $i\in D$, an initial value for the job arrival rate $\lambda'_i$ is sampled from a continuous uniform distribution between $0.1\mu_i$ and $\mu_i$.
        \item For each demand point $i\in D$, the actual traffic intensity $\rho_i$ is obtained by re-scaling the initial traffic intensity $\lambda'_i/\mu_i$, as follows:
        $$\rho_i:=\frac{\lambda'_i/\mu_i}{\sum_{i\in D}\lambda'_i/\mu_i}\rho.$$
        This ensures that $\sum_{i\in D}\rho_i=\rho$.
        \item For each demand point $i\in D$, the actual job arrival rate $\lambda_i$ is obtained as follows:
        $$\lambda_i:=\rho_i\mu_i.$$
        \item All of the $\lambda_i$, $\mu_i$ and $\tau$ values are re-scaled in order to ensure that $\sum_{i\in D}\lambda_i+\max\{\mu_1,...,\mu_d,\tau\}=1$ (as assumed in Section \ref{sec:Problem Formulation}) and then rounded to 2 significant figures. 
\end{itemize}
        
    \item For each demand point $i\in D$, the holding cost $c_i$ is sampled from a continuous uniform distribution between 0.1 and 0.9.
    \item In order to generate the switching rate $\tau$, we first define $\eta:=\tau/(\sum_{i\in D}\lambda_i)$ and generate the value of $\eta$ as follows:
        \begin{itemize}
        \item Sample a value $p$ from a continuous uniform distribution between 0 and 1.
        \item If $p < 0.5$, sample $\eta$ from a continuous uniform distribution between 0.1 and 1.
        \item If $p\geq 0.5$, sample $\eta$ from a continuous uniform distribution between 1 and 10.
        \end{itemize}
    We then define  $\tau := \eta \sum_{i\in D}\lambda_i$.\\
    \end{itemize}

\section{Simulation methods for the numerical experiments in Section \ref{sec:numerical}}\label{AppE}

In each of the problem instances \black{considered in Sections \ref{sec41} and \ref{sec42}}, the performances of the $K$-stop and ($K$ from $L$)-stop heuristics are estimated by simulating the discrete-time evolution of the uniformized MDP described in Section \ref{sec:Problem Formulation}. We also use a `common random numbers' method to ensure that job arrivals occur at the same times under each of these policies. More specifically, the simulation steps are as follows:
\begin{enumerate}
\item Generate a set of random system parameters as described in Appendix \ref{AppD}.
\item Set $N_0:=10,000$ as the length of the warm-up period and $N_1:=1,000,000$ as the length of the simulation.
\item Generate a list $Z$ of length $N_0+N_1=1,010,000$, consisting of uniformly-distributed random numbers between 0 and 1.
\item Consider each of the $K$-stop and ($K$ from $L$)-stop heuristics in turn. For each one, set $\vec{x}_0:=(1,(0,0,...,0))$ as the initial state and use the first $N_0$ random numbers in $Z$ to simulate events in the first $N_0$ time steps (this is done by sampling from the transition probability distribution described in (\ref{trans_prob_eqn})). This is the `warm-up period'. Let $\vec{y}$ denote the state reached at the end of the warm-up period. Then, use the remaining $N_1$ random numbers in $Z$ to simulate events during the next $N_1$ time steps, with the system beginning in state $\vec{y}$, and use the statistics collected during these time steps to quantify the heuristic's performance.\\
\end{enumerate}
The method for simulating the DVO heuristic is different, because (as discussed in Section \ref{sec:Index heuristic}) this heuristic is not directly compatible with the MDP formulation in the paper. To simulate the DVO heuristic, we simulate in continuous time as follows:
\begin{enumerate}
\item Use the same set of randomly-generated system parameters used for the other heuristics. 
\item Let $T_0=10,000$ and $T_1=1,000,000$ as the warm-up period length and the main simulation length, respectively.
\item Set $\vec{x}_0:=(1,(0,0,...,0))$ as the initial state. Note that, under the DVO heuristic, the set of decision epochs is not the same as under the other heuristics (see Section \ref{DVO_heuristic}). At each decision epoch we make a decision using the rules of the DVO heuristic and then simulate the time until the next decision epoch, which requires simulating from an exponential distribution (if the server processes a job or idles at an empty demand point) or an Erlang distribution (if the server switches to another demand point). Costs are accumulated based on the job counts $x_i$ at the various demand points in between decision epochs. This process continues until the total time elapsed exceeds $T_0$, at which point the warm-up period ends. Let $\vec{y}$ denote the state reached at the end of the warm-up period. We then carry out the main simulation in the same way, starting from state $\vec{y}$ and continuing for a further $T_1$ time units, at which point the process ends and the performance of the DVO heuristic is estimated by dividing the total costs incurred during the main simulation by the total time elapsed.\\
\end{enumerate}

\vspace{10pt}
\section{Method for testing the feasibility of dynamic programming (DP) for the numerical experiments in Section \ref{sec:numerical}} \label{AppF}
As explained in Section \ref{sec:numerical}, we aim to use DP (specifically, relative value iteration) to compute the optimal average cost $g^*$ whenever it is computationally feasible to do so. For a particular problem instance, we carry out the following steps in order to classify it as either `feasible' or `infeasible' and (for the feasible instances only) estimate the optimal average cost:
\begin{enumerate}
    \item If $d \geq 4$, classify this instance as infeasible.
    \item Otherwise (if $d \leq 3$), carry out the following steps:
        \begin{enumerate}
            \item Set $m=10$, $t = 0$ and $g^* = 0$.
            \item Consider a finite-state MDP in which the number of jobs present at any demand point is not allowed to exceed $m$, as described in the proof of Theorem \ref{stability_theorem}.  Let $M$ denote the size of the state space, given by
            $$M:=(d+n)(m+1)^d.$$
            (Recall that $d$ and $n$ are the numbers of demand points and intermediate stages in the network, respectively.) If $M\geq 1,000,000$, go to step (d). Otherwise, solve the finite-state MDP using DP, let $t$ denote the time taken (in seconds) and let $g^*$ denote the optimal average cost.
            \item If $t < 600$, increase $m$ by 10 and return to step (b). Otherwise, continue to step (d).
            \item If either (i) $g^*=0$, or (ii) the latest value of $g^*$ exceeds the previous value by more than $\epsilon$ (where we set $\epsilon=0.001$), then classify this instance as infeasible. Otherwise, classify it as feasible and let the latest value of $g^*$ be an approximation for the optimal average cost in the infinite-state MDP.
        \end{enumerate}
\end{enumerate}

\end{appendices}

\end{document}